\documentclass{amsart}
\usepackage{amssymb} 
\usepackage{color}
\usepackage[all]{xy}
\usepackage{enumerate}

 \newtheorem{theorem}{Theorem}[section]

 \newtheorem{proposition}[theorem]{Proposition}
 \newtheorem{lemma}[theorem]{Lemma}
 \newtheorem{corollary}[theorem]{Corollary}
 \newtheorem{remark}[theorem]{Remark}
 \newtheorem{definition}[theorem]{Definition}
 \newtheorem{notation}[theorem]{Notation}
 \newtheorem{example}[theorem]{Example}

 \newtheorem{claim}[theorem]{Claim}

 \numberwithin{equation}{section}

\def\cat{{\rm Cat}_{(\infty,1)}}

\def\op#1{{\rm Op}_{\infty/#1}}

\def\map#1{{\rm Map}_{#1}}

\def\bicat{{\rm Cat}_{(\infty,2)}}

\def\arrow{{\rm Ar}}
\def\oplaxarrow#1{{\rm Ar}^{\rm opl}(#1)}
\def\laxarrow#1{{\rm Ar}^{\rm lax}(#1)}

\def\oplaxarrowLR#1{{\rm Ar}^{\rm opl}_{L,R}(#1)}

\def\oplaxarrowsl#1{{\rm Ar}^{\rm opl}_{sL}(#1)}
\def\laxarrowsr#1{{\rm Ar}^{\rm lax}_{sR}(#1)}
\def\overlineoplaxarrow#1{\overline{\rm Ar}^{\rm opl}_{L,R}(#1)}

\def\funopl{{\rm Fun}^{\rm opl}}
\def\scriptx{\mathcal{X}}
\def\scripty{\mathcal{Y}}
\def\scriptz{\mathcal{Z}}
\def\scriptw{\mathcal{W}}
\def\scripto{\mathcal{O}}
\def\scriptp{\mathcal{P}}
\def\scriptc{\mathcal{C}}

\def\mapsc{{\rm Map}^{\rm sc}}
\def\ts{{\rm ts}}
\def\setsc{{\rm Set}_{\Delta}^{\rm sc}}
\def\tw{{\rm Tw}}
\def\twl#1{{\rm Tw}^l(#1)}
\def\twlsr#1{{\rm Tw}^l_{sR}(#1)}

\def\twr#1{{\rm Tw}^r(#1)}
\def\twrsl#1{{\rm Tw}^r_{sL}(#1)}
\def\mapsc{{\rm Map}^{\rm sc}}
\def\ts{{\rm ts}}
\def\setsc{{\rm Set}_{\Delta}^{\rm sc}}
\def\tw{{\rm Tw}}

\def\setjoy{{\rm Set}_{\Delta}^{\rm Joy}}

\def\setkan{{\rm Set}_{\Delta}^{\rm Kan}}
\def\setdelta{{\rm Set}_{\Delta}}

\def\tsbullet{{\rm ts}^{\bullet}}

\def\tsscbullet{{\rm ts}^{\bullet}}
\def\tssc{{\rm ts}}
\def\tsscn{{\rm ts}^{n}}
\def\tssczero{{\rm ts}^{0}}
\def\tsscone{{\rm ts}^{1}}

\def\wtsscone{{\rm ts}^{1}_{\diamondsuit}}
\def\tsscplus{{\rm ts}^{\bullet}_{\rm UR}}
\def\tsscplusno{{\rm ts}_{\rm UR}}
\def\tsscnplus{{\rm ts}^n_{\rm UR}}
\def\tssconeplus{{\rm ts}^{1}_{\rm UR}}
\def\tssconeminus{{\rm ts}^{1}_{\rm LL}}
\def\tsscminus{{\rm ts}^{\bullet}_{\rm LL}}
\def\tsscminusno{{\rm ts}_{\rm LL}}
\def\tsscnminus{{\rm ts}^{n}_{\rm LL}}
\def\TSop{{\rm TS}}

\def\TSopl{{\rm TS}_{L}}
\def\TSoprl{{\rm TS}_{L,R}}
\def\TSoprlcocart{{\rm TS}_{L,R}^{\rm cocart}}
\def\TSoprlpair{{\rm TS}_{L,R}^{\rm pair}}

\def\twoop{\mbox{\rm\scriptsize $2$-op}}
\def\bd{{\rm Bd}}

\begin{document}

\title{
Map monoidales
and duoidal $\infty$-categories}
\author{Takeshi Torii}
\address{Department of Mathematics, 
Okayama University,
Okayama 700--8530, Japan}
\email{torii@math.okayama-u.ac.jp}


\date{June 1, 2024 (version~1.0)}

\begin{abstract}
In this paper
we give an example of duoidal $\infty$-categories.
We introduce map $\mathcal{O}$-monoidales
in an $\mathcal{O}$-monoidal $(\infty,2)$-category
for an $\infty$-operad $\mathcal{O}^{\otimes}$.
We show that the endomorphism mapping
$\infty$-category of a map $\mathcal{O}$-monoidale
is a coCartesian $(\Delta^{\rm op},\mathcal{O})$-duoidal
$\infty$-category.
After that,
we introduce a convolution product
on the mapping $\infty$-category
from an $\mathcal{O}$-comonoidale to
an $\mathcal{O}$-monoidale.
We show that the $\mathcal{O}$-monoidal
structure on the duoidal endomorphism mapping $\infty$-category
of a map $\mathcal{O}$-monoidale
is equivalent to the convolution product
on the mapping $\infty$-category 
from the dual $\mathcal{O}$-comonoidale
to the map $\mathcal{O}$-monoidale.
\end{abstract}

\maketitle

\section{Introduction}

Duoidal categories are categories equipped with
two monoidal products in which one is 
(op)lax monoidal with respect to the other.
Duoidal categories are generalizations of braided monoidal categories.
Braided monoidal categories can be regarded as pseudomonoids
in the monoidal $2$-category of monoidal categories
and strong monoidal functors.
On the other hand,
duoidal categories can be regarded as pseudomonoids
in the monoidal $2$-category of monoidal categories and
lax monoidal functors.
The notion of duoidal category was introduced
by Aguiar-Mahajan by name of $2$-monoidal category
in \cite{Aguiar-Mahajan}.
The terminology of duoidal category was proposed by
Street in \cite{Street}.

Examples and applications
of duoidal categories
are presented
in \cite{Aguiar-Mahajan}
and \cite{Street}.
We can make sense of monoidal categories enriched over
duoidal categories.
In \cite{Batanin-Markl1, Batanin-Markl2}
duoidal categories play fundamental role 
in their theory of centers and homotopy centers of
monoids in monoidal categories enriched
over duoidal categories,
and in developing of a generalization of Deligne's conjecture.
A general theory of duoidal categories
have further been developed in \cite{Booker-Street}.
We can consider bimonoids in duoidal categories,
which are generalizations of bialgebras.
There are several equivalent conditions
for bialgebras to be Hopf algebras. 
The fundamental theorem of Hopf modules 
for bimonoids in duoidal categories was studied
in \cite{BCZ}.
The notion of bimonoids in duoidal categories 
have been generalized to that of weak bimonoids
in \cite{Chen-Bohm}.
Hopf conditions 
have further been studied 
for bimonoids in duoidal endomorphism hom-categories
of map monoidales in monoidal bicategories
in \cite{Bohm-Lack}.

As explained in \cite{Bohm-Lack},
the endomorphism hom-category
of a map monoidale in a monoidal bicategory
admits the structure of a duoidal category.
One of the monoidal products
is given by the composition of morphisms
and the other is by convolution product.
In this paper we generalize
this construction in the setting of higher categories.

In higher category theory
we replace the category of sets in ordinary category theory
to the homotopy theory of spaces
based on the homotopy hypothesis \cite{Lurie1}.
In other words,
the theory of $(\infty,0)$-categories
is equivalent to the homotopy theory of spaces.
An $(\infty,n)$-category is an $\infty$-category
whose $k$-morphisms are invertible for $k>n$.
On way to view an $(\infty,n)$-category
is to regard it as an $\infty$-category
enriched over $(\infty,n-1)$-categories.
In particular,
an $(\infty,1)$-category is an $\infty$-category
enriched over spaces. 
Another way to view an $(\infty,n)$-category
is to regard it as an
$n$-category up to coherent homotopy.
By taking the set of components of mapping spaces,
we can associate to an $(\infty,1)$-category
its homotopy category which is an ordinary category.
An $(\infty,2)$-category
is regarded as an $\infty$-category
enriched over $(\infty,1)$-categories.
By taking homotopy categories
of mapping $(\infty,1)$-categories,
we obtain the homotopy bicategory of
an $(\infty,2)$-category.

In \cite{Torii1}
we introduced the notion of duoidal $\infty$-categories,
which are counterparts of duoidal categories
in the setting of higher categories.
In \cite{Torii2}
we generalized duoidal $\infty$-categories
to higher monoidal $\infty$-categories
and studied their duality property.
Higher monoidal categories
are generalizations of duoidal categories, 
which are also introduced in \cite{Aguiar-Mahajan}.
Higher monoidal categories have multiple monoidal products
together with multiple interchange laws.
Higher monoidal $\infty$-categories
are counterparts of higher monoidal categories
in higher category theory.
In \cite{Torii3}
we gave an example of duoidal $\infty$-categories
of operadic modules.
Lurie \cite{Lurie2} showed that the $(\infty,1)$-category
of $\scripto$-$A$-modules
in $\scriptc$ admits
the structure of an $\scripto$-monoidal $\infty$-category,
where $\scripto^{\otimes}$ is
a small coherent symmetric $\infty$-operad,
$\scriptc$ is an $\scripto$-monoidal $(\infty,1)$-category
with a sufficient supply of colimits,
and $A$ is an $\scripto$-algebra in $\scriptc$.
By generalizing this construction,
we showed that if $\mathcal{P}^{\otimes}$
is a symmetric $\infty$-operad,
$\scripto^{\otimes}$ is a small coherent
symmetric $\infty$-operad,
$\mathcal{C}$
is a $\scriptp\otimes\scripto$-monoidal
$\infty$-category
with a sufficient supply of colimits,
and $A$ is a $\scriptp\otimes\scripto$-algebra
in $\scriptc$,
then the $(\infty,1)$-category
of $\scripto$-$A$-modules in $\scriptc$
admits the structure
of a mixed $(\scriptp,\scripto)$-duoidal
$\infty$-category,
where $\scriptp^{\otimes}\otimes\scripto^{\otimes}$
is the Boardman-Vogt tensor product
of $\scriptp^{\otimes}$ and $\scripto^{\otimes}$.

The purpose of this paper is
to give another example of duoidal $\infty$-categories.
We show that we can construct a duoidal
$\infty$-category by the same procedure as
in \cite{Bohm-Lack}
in the setting of higher categories.
We consider an $\scripto$-monoidal
$(\infty,2)$-category $\scriptz$
for an $\infty$-operad $\scripto^{\otimes}$
over a perfect operator category.
We introduce an $\scripto$-monoidale
as an $\scripto$-monoid object
in the underlying $(\infty,1)$-category of $\scriptz$.
We can consider
the $1$-full sub $(\infty,2)$-category $\scriptz^L$
spanned by left adjoint $1$-morphisms,
which inherits $\scripto$-monoidal structure
from $\scriptz$.
A map $\scripto$-monoidale
is defined to be an $\scripto$-monoidale in $\scriptz^L$.

For an object $A$ in $\scriptz$,
the endomorphism mapping $\infty$-category $\scriptz(A,A)$
is a monoidal $(\infty,1)$-category
by composition of morphisms.
Let $\Delta^{\rm op}$ be the opposite category of
nonempty totally ordered sets.
We regard it as an $\infty$-operad over
the perfect operator category $\mathbf{O}$
of finite ordered sets.
We formulate the composition monoidal structure
on $\scriptz(A,A)$ by using
the $\infty$-operad $\Delta^{\rm op}$.

The first main theorem of this paper is as follows.

\begin{theorem}[{cf.~Theorem~\ref{thm:Z-A-A-duoidal-structure}}]
\label{theorem:first-main-theoorem}
Let $\scriptz$ be an $\scripto$-monoidal
$(\infty,2)$-category and let
$A$ be a map $\scripto$-monoidale in $\scriptz$.
The endomorphism mapping $\infty$-category $\scriptz(A,A)$
admits the structure of a coCartesian
$(\Delta^{\rm op}, \scripto)$-duoidal $\infty$-category,
\end{theorem}

Next,
we compare
the $\scripto$-monoidal structure
on $\scriptz(A,A)$ with convolution product.   
We call an $\scripto$-comonoid object
in the underlying $(\infty,1)$-category of $\scriptz$
an $\scripto$-comonoidale.
For an $\scripto$-monoidale $A$
and an $\scripto$-comonoidale $C$,
we can define a convolution product
on the mapping $\infty$-category $\scriptz(C,A)$.
Informally speaking,
when $\scripto^{\otimes}$ is the associative $\infty$-operad,
the convolution product is given by
\[ C\stackrel{c}{\longrightarrow}
   C\otimes C\stackrel{f\otimes g}{\longrightarrow}
   A\otimes A\stackrel{m}{\longrightarrow} A\]
for $f,g\in \scriptz(C,A)$,
where $m$ is the product on $A$ and $c$ is the coproduct
on $C$.
The convolution product
gives $\scriptz(C,A)$
the structure of an $\scripto$-monoidal
$(\infty,1)$-category.

We denote by $\scriptz^R$
the $1$-full sub $(\infty,2)$-category of $\scriptz$
spanned by right adjoint $1$-morphisms.
By taking adjoint morphisms,
there is a dual equivalence 
between the $(\infty,1)$-category of
$\scripto$-monoidales in $\scriptz^L$ 
and the $(\infty,1)$-category
of $\scripto$-comonoidales
in $\scriptz^R$.
We denote by $A^*$
the $\scripto$-comonoidale in $\scriptz^R$    
corresponding to a map $\scripto$-monoidale $A$
under this equivalence.
Since the underlying object of $A^*$
is equivalent to that of $A$,
the underlying $(\infty,1)$-category
of $\scriptz(A^*,A)$ is equivalent
to that of $\scriptz(A,A)$.

The second main theorem of this paper is as follows.

\begin{theorem}
Let $\scriptz$ be an $\scripto$-monoidal $(\infty,2)$-category
and let $A$ be a map $\scripto$-monoidale 
in $\scriptz$.
The $\scripto$-monoidal structure
on the coCartesian $(\Delta^{\rm op},\scripto)$-duoidal
$\infty$-category $\scriptz(A,A)$
is equivalent to the convolution product
on $\scriptz(A^*,A)$.
\end{theorem}


The organization of this paper is as follows:
The paper is divided into two parts.
In Part 1 we introduce a notion of map monoidales
and construct a duoidal structure on
the endomorphism mapping $\infty$-category of a map monoidale.
In \S\ref{section:review-duoidal-categoreis}
we review duoidal categories
in the setting of higher categories.
We recall the definition of
(virtual) coCartesian $(\scripto,\scriptp)$-duoidal $\infty$-categories,
where $\scripto^{\otimes}$ and $\scriptp^{\otimes}$
are $\infty$-operads over perfect operator categories.
In \S\ref{section:scaled-simplicial-set}
we review scaled simplicial sets
which are models of $(\infty,2)$-categories.
We recall the model structure on
the category of scaled simplicial sets
whose underlying $\infty$-category
is equivalent to the theory of
$(\infty,2)$-categories.
In \S\ref{section:adjunctions-infty-bicategory}
we study adjoint morphisms in $(\infty,2)$-categories.
We construct $1$-full sub $(\infty,2)$-categories
spanned by left and right adjoint $1$-morphisms,
respectively.
We show that there is a natural dual equivalence
between the underlying $(\infty,1)$-categories of them.
In \S\ref{section:map-monoidale}
we introduce a map $\scripto$-monoidale
in an $\scripto$-monoidal $(\infty,2)$-category,
where $\scripto^{\otimes}$ is an $\infty$-operad 
over a perfect operator category.
We define a map $\scripto$-monoidale
to be an $\scripto$-monoidale in $\scriptz^L$.
We show that there is a dual equivalence
between the $(\infty,1)$-category of $\scripto$-monoidales
in $\scriptz^L$ and
the $(\infty,1)$-category
of $\scripto$-comonoidales in $\scriptz^R$ 
by taking adjoint morphisms.
In \S\ref{section:map-duoidal}
we prove Theorem~\ref{theorem:first-main-theoorem}.
First,
we show that the endomorphism mapping $\infty$-category
$\scriptz(A,A)$
of an $\scripto$-monoidale $A$ is a virtual coCartesian
$(\Delta^{\rm op},\scripto)$-duoidal $\infty$-category.
After that,
we show that it is in fact a coCartesian
$(\Delta^{\rm op},\scripto)$-duoidal $\infty$-category
when $A$ is a map $\scripto$-monoidale.

In Part 2
we compare the $\scripto$-monoidal structure
on $\scriptz(A,A)$ with convolution product.
In \S\ref{section:convolution_product}
we introduce a convolution product
on the mapping $\infty$-category.
For an $\scripto$-monoidale $A$
and an $\scripto$-comonoidale $C$,
we show that the mapping $\infty$-category
$\scriptz(C,A)$ has the structure
of an $\scripto$-monoidal $(\infty,1)$-category
by convolution product.
In order to compare the $\scripto$-monoidal structure
on $\scriptz(A,A)$ with convolution product
when $A$ is a map $\scripto$-monoidale,
we need to establish an equivalence
between a wide subcategory
$\twlsr{\scriptx}$
of the twisted arrow $\infty$-category 
and a wide subcategory $\oplaxarrowsl{\scriptx}$
of the oplax arrow $\infty$-category
of an $(\infty,2)$-category $\scriptx$.
For this purpose,
in \S\ref{section:twisted-squares}
we introduce a twisted square $\infty$-category
$\TSop(\scriptx)$,
which is a generalization of
twisted arrow $\infty$-category.
For a fibrant scaled simplicial set $X$
which represents $\scriptx$,
we construct a simplicial space
$\TSop(X)_{\bullet}$
which represents $\TSop(\scriptx)$.
We defer to \S\ref{section:proof-twisted-sq-complete-Segal}
a proof of 
$\TSop(X)_{\bullet}$ being a complete Segal space.
In \S\ref{section:twl-arrow-equivalence}
we show that there is an equivalence 
of $(\infty,1)$-categories
between $\twlsr{\scriptx}$
and $\oplaxarrowsl{\scriptx}$.
For this purpose,
we construct a perfect pairing between them
by taking a suitable subcategory of
$\TSop(\scriptx)$.
Finally,
in \S\ref{section:duoidal-convolution}
we show that the $\scripto$-monoidal
structure on $\scriptz(A,A)$
is equivalent to the convolution product on
$\scriptz(A^*,A)$,
where $A^*$ is the dual $\scripto$-comonoidale
associated to the map $\scripto$-monoidale $A$.


\begin{notation}\rm
We denote by
$\setdelta$
the category of simplicial sets.
We write $\setkan$ and $\setjoy$
for the category $\setdelta$
equipped with the Kan and Joyal model structures,
respectively.
In this paper
we mean that an $\infty$-groupoid is a Kan complex,
that an $\infty$-category is a quasi-category,
and that an $\infty$-bicategory is
a fibrant scaled simplicial set.
We let ${\rm rev}: \Delta^{\rm op}\to\Delta^{\rm op}$
be a functor given by assigning to a nonempty finite ordered
set its reverse ordered set. 
For a simplicial object $S_{\bullet}$,
we denote by $S_{\bullet}^{\rm rev}$
the composite functor $S_{\bullet}\circ {\rm rev}$.

We denote by $\mathcal{S}$
the $\infty$-category of (small) spaces,
by $\cat$ the $\infty$-category
of (small) $(\infty,1)$-categories,
and
by $\bicat$ the $\infty$-category
of (small) $(\infty,2)$-categories.
We denote by
$(-)^{\simeq}: \cat\to\mathcal{S}$
the right adjoint
to the inclusion functor
$\mathcal{S}\to\cat$,
and by
$u_1(-): \bicat\to\cat$
the right adjoint to the inclusion functor
$\cat\to\bicat$.
For an $\infty$-category $\mathcal{C}$
which admits finite products,
we denote by $\mathcal{C}^{\times}$
the associated Cartesian symmetric monoidal $\infty$-category.

For an $(\infty,1)$-category $\mathcal{C}$
with objects $c$ and $d$,
we denote by $\map{\mathcal{C}}(c,d)$
the mapping space from $c$ to $d$.
For an $(\infty,2)$-category $\mathcal{X}$ 
with objects $x$ and $y$,
we denote by $\mathcal{X}(x,y)$
the mapping $(\infty,1)$-category from $x$ to $y$.

For a nonempty finite totally ordered set $S$,
we denote by $\Delta^S$ the $(|S|-1)$-dimensional simplex
with $S$ as the set of vertices.
For $s\in S$,
$\Delta^{S-{\{s\}}}$ 
is the codimension $1$ face of $\Delta^S$
opposite to the vertex $s$.
For a subset $\emptyset \neq N\varsubsetneqq S$,
we set
\[ \Lambda^S_N=\bigcup_{s\in S-N}\Delta^{S-\{s\}}. \]
For $[r]=\{0<1<\cdots<r\}$,
we write $\Delta^r$ and $\Lambda^r_i$ for
$\Delta^{[r]}$ and $\Lambda^{[r]}_{\{i\}}$,
respectively.
\end{notation}

\part{Duoidal endomorphism mapping $\infty$-categories
of map monoidales}
\label{part:endomorphism-duoidal}

In ordinary category theory
a map monoidale is a pseudomonoid
in a monoidal bicategory whose
unit and multiplication
morphisms are left adjoint.
The endomorphism category of a map monoidale
admits the structure of a duoidal category.
In Part~\ref{part:endomorphism-duoidal}
we consider counterparts
of this construction in higher category theory.

In \S\ref{section:map-monoidale}
we introduce a notion of a map $\scripto$-monoidale in
an $\scripto$-monoidal $(\infty,2)$-category $\scriptz$,
where $\scripto$ is an $\infty$-operad over a
perfect operator category.
In \S\ref{section:map-duoidal}
we show that the endomorphism mapping $\infty$-category
of a map $\scripto$-monoidale in $\scriptz$
admits the structure of
a coCartesian $(\Delta^{\rm op},\scripto)$-duoidal
$\infty$-category.

For this purpose,
in \S\ref{section:review-duoidal-categoreis}
we recall the notion of duoidal $\infty$-categories
introduced in \cite{Torii1, Torii2}.
In \S\ref{section:scaled-simplicial-set}
we recall scaled simplicial sets
which are models of $(\infty,2)$-categories.
In \S\ref{section:adjunctions-infty-bicategory}
we study adjoint morphisms in $(\infty,2)$-categories.

\section{Duoidal $\infty$-categories}
\label{section:review-duoidal-categoreis}

In ordinary category theory
a duoidal category has two monoidal products
in which one is (op)lax monoidal with respect
to the other.
In this section
we review the notion of duoidal category
in the setting of higher category theory
(cf.~\cite{Torii1, Torii2}).

Let $\Phi$ be a prefect operator category 
in the sense of \cite[Definitions~1.2 and 4.6]{Barwick}.
Associated to $\Phi$,
we have the Leinster category $\Lambda(\Phi)$
equipped with collections of inert morphisms and active morphisms.
An $\infty$-operad over $\Phi$
is a functor $p: \mathcal{O}^{\otimes}\to \Lambda(\Phi)$
of $\infty$-categories satisfying certain conditions
(see \cite[Definition~7.8]{Barwick}
for precise definition).
By definition,
a perfect operator category $\Phi$
has a final object $*$.
We denote by $\mathcal{O}$
the fiber $\mathcal{O}^{\otimes}_*$ at $*\in \Lambda(\Phi)$
and say that it is the underlying $\infty$-category of 
$\mathcal{O}^{\otimes}$. 
A morphism of $\mathcal{O}^{\otimes}$
is said to be inert if it is a $p$-coCartesian morphism
over an inert morphism,
and active if it covers an active morphism
of $\Lambda(\Phi)$.
A morphism of $\infty$-operads over $\Phi$
between $\mathcal{O}^{\otimes}\to\Lambda(\Phi)$
and $\mathcal{P}^{\otimes}\to\Lambda(\Phi)$
is a functor $f: \mathcal{O}^{\otimes}\to \mathcal{P}^{\otimes}$
over $\Lambda(\Phi)$ that preserves inert morphisms.
We denote by
$\op{\Lambda(\Phi)}$
the $\infty$-category of $\infty$-operads
over $\Phi$.
We notice that
$\op{\Lambda(\Phi)}$ admits finite products.

\begin{example}
\rm
Let $\mathbf{F}$ be the category of finite sets.
By \cite[Example~4.9.3]{Barwick},
$\mathbf{F}$ is a perfect operator category.
By \cite[Example~6.5]{Barwick},
$\Lambda(\mathbf{F})$ 
is equivalent to the category ${\rm Fin}_*$ 
of pointed finite sets.
By \cite[Example~7.9]{Barwick},
the notion of $\infty$-operads over $\mathbf{F}$
coincides with that of Lurie's 
$\infty$-operads
in \cite[Chapter~2]{Lurie2}. 
\end{example}

\begin{example}
\rm
Let $\mathbf{O}$ be the category of ordered finite sets.
By \cite[Example~4.9.2]{Barwick},
$\mathbf{O}$ is a perfect operator category.
By \cite[Example~6.6]{Barwick},
$\Lambda(\mathbf{O})$ is equivalent to $\Delta^{\rm op}$.
The notion of $\infty$-operads over $\mathbf{O}$
coincides with that of non-symmetric $\infty$-operads
in \cite[\S3]{Gepner-Haugseng}. 
\end{example}

Let $C^{\otimes}\to \mathcal{O}^{\otimes}$
be a map of $\infty$-operads over $\Phi$.
We say that it is an $\mathcal{O}$-monoidal $\infty$-category if 
it is a coCartesian fibration.
In this case
we have a multiplication map
\[ \otimes_{\phi}:
   \prod_{i\in |I|}C_{a_i}^{\otimes}
   \stackrel{\simeq}{\longleftarrow}
   C^{\otimes}_a
   \stackrel{\phi_*}{\longrightarrow}
   C^{\otimes}_b\]
for an active morphism $\phi: a\to b$ in $\mathcal{O}^{\otimes}$,
where $p(a)=I$, $a\simeq (a_i)_{i\in |I|}$,
and $a_i,b\in \mathcal{O}$.

Let $C$ and $D$
be $\mathcal{O}$-monoidal $\infty$-categories,
and let 
$g: C^{\otimes}\to D^{\otimes}$
be a map of $\infty$-categories over $\mathcal{O}^{\otimes}$.
We say that
$g$ is a lax $\mathcal{O}$-monoidal functor
if it is a map of $\infty$-operads over $\mathcal{O}^{\otimes}$,
that is,
it preserves inert morphisms.
We say that
$g$ is a strong $\mathcal{O}$-monoidal functor
if it preserves coCartesian morphisms.
We denote by
$\mathsf{Mon}_{\mathcal{O}}^{\rm lax}(\cat)$
the subcategory of ${\rm Cat}_{(\infty,1)/\mathcal{O}^{\otimes}}$
spanned by $\mathcal{O}$-monoidal $\infty$-category
and lax $\mathcal{O}$-monoidal functors.
We denote by $\mathsf{Mon}_{\mathcal{O}}(\cat)$
the wide subcategory of $\mathsf{Mon}_{\mathcal{O}}^{\rm lax}(\cat)$
spanned by strong $\mathcal{O}$-monoidal functors.
We notice that 
$\mathsf{Mon}_{\mathcal{O}}^{\rm lax}(\cat)$
and 
$\mathsf{Mon}_{\mathcal{O}}(\cat)$
admit finite products.

As in \cite[\S2.4.2]{Lurie2},
we can identify an $\mathcal{O}$-monoidal
$\infty$-category with a functor
from $\mathcal{O}^{\otimes}$ to $\cat$.
Let $Y$ be an $\infty$-category
with finite products,
and let $M$ be a functor $\mathcal{O}^{\otimes}\to Y$
of $\infty$-categories.
For any $a\in\mathcal{O}^{\otimes}$
with $p(a)=I$,
we can take a family of $p$-coCartesian
morphisms $\{a\to a_i|\ i\in |I|\}$
lying over the inert morphisms
$\{\rho_i: I\to \{i\}|\ i\in |I|\}$.
The family $\{a\to a_i|\ i\in |I|\}$
induces a morphism $M(a)\to \prod_{i\in |I|}M(a_i)$
in $Y$,
which we call a Segal map.
We say that $M$ is an $\mathcal{O}$-monoid object in $Y$
if the Segal map is an equivalence
for each $a\in\mathcal{O}^{\otimes}$.
We denote by
${\rm Mon}_{\mathcal{O}}(Y)$
the full subcategory of ${\rm Fun}(\mathcal{O}^{\otimes},Y)$
spanned by $\mathcal{O}$-monoid objects.
There is an equivalence
of $\infty$-categories between 
$\mathsf{Mon}_{\mathcal{O}}(\cat)$
and
${\rm Mon}_{\mathcal{O}}(\cat)$
(cf.~\cite[Remark~2.4.2.6]{Lurie2}).

Now,
we define a notion of virtual
coCartesian duoidal $\infty$-category.
Let $p_{\mathcal{O}}: \mathcal{O}^{\otimes}\to \Lambda(\Phi)$ 
and $p_{\mathcal{P}}: \mathcal{P}^{\otimes}\to \Lambda(\Psi)$
be $\infty$-operads over perfect
operator categories $\Phi$ and $\Psi$, respectively.
We denote by
$\pi_{\mathcal{O}}$ and $\pi_{\mathcal{P}}$
the projections from
$\mathcal{O}^{\otimes}\times\mathcal{P}^{\otimes}$
to the factors, respectively.

\begin{definition}
[Definition of virtual coCartesian Duoidal $\infty$-categories]
\rm
Let 
\[ f: X \to \mathcal{O}^{\otimes}\times\mathcal{P}^{\otimes}. \]
be a categorical fibration of $\infty$-categories.
We say that $f$
is a virtual coCartesian $(\mathcal{O},\mathcal{P})$-duoidal
$\infty$-category 
if it satisfies the following conditions:

\begin{enumerate}
\item The composite $\pi_{\mathcal{O}}\circ f: 
X\to \mathcal{O}^{\otimes}$
is a coCartesian fibration,
and $f$ carries $(\pi_{\mathcal{O}}\circ f)$-coCartesian 
morphisms to $\pi_{\mathcal{O}}$-coCartesian morphisms.
\item
For each $a\in \mathcal{O}^{\otimes}$,
the restriction
$p_{a}: X_{a}\to \mathcal{P}^{\otimes}$
is a morphism of $\infty$-operads
over $\Lambda(\Psi)$.

\item
For each morphism
$a\to a'$ in $\mathcal{O}^{\otimes}$,
the induced map
$X_{a}\to X_{a'}$
over $\mathcal{P}^{\otimes}$
is a morphism of $\infty$-operads
over $\mathcal{P}^{\otimes}$.
\item
For each $a\in \mathcal{O}^{\otimes}$
with $p_{\mathcal{O}}(a)=I$,
the Segal morphism
\[ X_{a}\longrightarrow
       \prod_{i\in |I|}{}^{\mathcal{P}^{\otimes}} X_{a_i} \]
is an equivalence in the $\infty$-category
$\op{\mathcal{P}^{\otimes}}$
of $\infty$-operads over $\mathcal{P}^{\otimes}$,
where 
$\{a\to a_i|\ i\in |I|\}$
is a family of $p_{\mathcal{O}}$-coCartesian morphisms
lying over the inert morphisms $\{\rho_i: I\to\{i\}|\ i\in |I|\}$,
and the right hand side is a product 
in $\op{\mathcal{P}^{\otimes}}$.
\end{enumerate}
\end{definition}

\begin{remark}\label{remark:cor-virtual-duoidal-category}
\rm
A virtual coCartesian
$(\mathcal{O},\mathcal{P})$-duoidal $\infty$-category 
corresponds to an $\mathcal{O}$-monoid object
of the Cartesian symmetric monoidal
$\infty$-category ${\rm Op}_{\infty/\mathcal{P}^{\otimes}}$
under the straightening along
the map $\pi_{\mathcal{O}}\circ f: X\to \mathcal{O}^{\otimes}$.
\end{remark}

Next,
we recall the definition
of coCartesian duoidal $\infty$-categories.

\begin{definition}
[Definition of coCartesian Duoidal $\infty$-categories
({cf.~\cite[Definition~4.23]{Torii1}
and \cite[Definition~3.14]{Torii2}})]
\rm
Let 
\[ f: X \to \mathcal{O}^{\otimes}\times\mathcal{P}^{\otimes}. \]
be a categorical fibration of $\infty$-categories. 
We say that $f$
is a coCartesian $(\mathcal{O},\mathcal{P})$-duoidal
$\infty$-category 
if it satisfies the following conditions:

\begin{enumerate}

\item
The map $f$ is a virtual coCartesian $(\mathcal{O},\mathcal{P})$-duoidal
$\infty$-category. 

\item
For each $a\in \mathcal{O}^{\otimes}$,
the restriction
$p_{a}: X_{a}\to \mathcal{P}^{\otimes}$
is a $\mathcal{P}$-monoidal $\infty$-category.


\end{enumerate}
\end{definition}

\begin{remark}\label{remark:cor-duoidal-category}
\rm
A coCartesian $(\mathcal{O},\mathcal{P})$-duoidal $\infty$-category 
corresponds to an $\mathcal{O}$-monoid object
of the Cartesian symmetric monoidal
$\infty$-category $\mathsf{Mon}_{\mathcal{P}}^{\rm lax}(\cat)$
under the straightening along
the map $\pi_{\mathcal{O}}\circ f: X\to \mathcal{O}^{\otimes}$.
\end{remark}

Informally speaking,
a coCartesian $(\mathcal{O},\mathcal{P})$-duoidal
$\infty$-category has $\mathcal{O}$-monoidal 
and $\mathcal{P}$-monoidal products in which
the $\mathcal{O}$-monoidal structure is lax monoidal
with respect to the $\mathcal{P}$-monoidal structure
and
the $\mathcal{P}$-monoidal structure
is oplax monoidal with respect to
the $\mathcal{O}$-monoidal structure.

\section{Scaled simplicial sets and $(\infty,2)$-categories}
\label{section:scaled-simplicial-set}

In this section
we review scaled simplicial sets
which are models of $(\infty,2)$-categories.
We recall the model structure on
the category of scaled simplicial sets introduced in \cite{Lurie3} 
whose underlying $\infty$-category
is equivalent to the $\infty$-category $\bicat$
of $(\infty,2)$-categories.

First,
we recall the definition of scaled simplicial sets.

\begin{definition}[{\cite[Definition~3.1.1]{Lurie3}})]\rm
A scaled simplicial set is a pair $(\overline{X},T_X)$
in which $\overline{X}$ is a simplicial set
and $T_X$ is a set of $2$-simplices of $\overline{X}$
that contains all the degenerate $2$-simplices.
We say that a $2$-simplex is thin
if it is contained in $T_X$.
For scaled simplicial sets
$X=(\overline{X},T_X)$ and $Y=(\overline{Y},T_Y)$, 
a morphism $f: X\to Y$ of scaled simplicial sets
is a map
$\overline{f}: \overline{X}\to \overline{Y}$
of simplicial sets that
carries $T_X$ into $T_Y$.
We denote by ${\rm Set}_{\Delta}^{\rm sc}$
the category of scaled simplicial sets.
\end{definition}

For a simplicial set $S$,
we have two canonical scaled simplicial sets $S_{\sharp}$
and $S_{\flat}$.
The scaled simplicial set $S_{\sharp}$
has the underlying simplicial set $S$
equipped with all $2$-simplices as thin triangles.
On the other hand,
$S_{\flat}$ has the underlying simplicial set $S$
equipped with all degenerate $2$-simplices as thin triangles.

Now, we recall scaled anodyne maps of
scaled simplicial sets which characterize fibrant objects
in $\setsc$.

\begin{definition}[{\cite[Definition~3.1.3]{Lurie3}}]\rm
The collection of scaled anodyne maps is
the weakly saturated class of morphisms
of ${\rm Set}_{\Delta}^{\rm sc}$ generated
by the following maps:
\begin{enumerate}

\item[(An1)]
For each $0<i<n$,
the inclusion map
\[ (\Lambda^n_i,\{\Delta^{\{i-1,i,i+1\}}\}|_{\Lambda^n_i}\cup
   \{{\rm degenerate}\})
   \longrightarrow
   (\Delta^n,\{\Delta^{\{i-1,i,i+1\}}\}\cup
    \{{\rm degenerate}\}).\]

\item[(An2)]
The inclusion map
\[ (\Delta^4,T)\longrightarrow
   (\Delta^4,T\cup \{\Delta^{\{0,3,4\}},\Delta^{0,1,4}\}),\]
where
\[ T=\{\Delta^{\{0,2,4\}},\Delta^{\{1,2,3\}},
       \Delta^{\{0,1,3\}},\Delta^{\{1,3,4\}},
       \Delta^{\{0,1,2\}}\}\cup\{{\rm degenerate}\}.\]

\item[(An3)] 
For $n>2$,
the inclusion map
\[ (\Lambda^n_0\coprod_{\Delta^{\{0,1\}}}\Delta^0,
   \{\Delta^{\{0,1,n\}}\}\cup\{{\rm degenerate}\})\longrightarrow
   (\Delta^n\coprod_{\Delta^{\{0,1\}}}\Delta^0,
   \{\Delta^{\{0,1,n\}}\}\cup\{{\rm degenerate}\} )   .\]
  
\end{enumerate}
\end{definition}

A model structure
on ${\rm Set}_{\Delta}^{\rm sc}$
was introduced in \cite{Lurie3} 
and showed that the underlying $\infty$-category
of ${\rm Set}_{\Delta}^{\rm sc}$ is
equivalent to the $\infty$-category
$\bicat$ of $(\infty,2)$-categories. 

\begin{theorem}[{cf.~\cite[Theorem~4.2.7]{Lurie3},
\cite[Theorem~1.29]{GHL},
and \cite[Theorem~1.12]{GHL2}}]
There exists a model structure on ${\rm Set}_{\Delta}^{\rm sc}$
whose cofibrations are the monomorphisms and
whose fibrant objects are those objects
that have the right lifting properties with respect
to the scaled anodyne maps.
\end{theorem}

\begin{definition}\rm
We call a fibrant scaled simplicial set
an $\infty$-bicategory. 
\end{definition}

\begin{remark}\rm
In \cite[Definition~4.1.1]{Lurie3}
scaled simplicial sets satisfying
the right lifting property
with respect to the scaled anodyne maps
are referred to as weak $\infty$-bicategories,
and fibrant objects in $\setsc$
are referred to as $\infty$-bicategories.
However, it is shown that these two notations coincide in
\cite[Theorem 5.1]{GHL}.
\end{remark}

We recall that 
$\setjoy$
is the category of simplicial sets
equipped with the Joyal model structure.
By assigning the scaled simplicial set $S_{\sharp}$ to
a simplicial set $S$,
we obtain a left Quillen functor $(-)_{\sharp}: \setjoy \to \setsc$.
We denote by
\[ u_1: \setsc \longrightarrow \setjoy  \]
the right adjoint to $(-)_{\sharp}$,
which assigns to a scaled simplicial set $X=(\overline{X},T_X)$
the subcomplex of $\overline{X}$ spanned by those simplices
whose $2$-dimensional faces are all thin.
Since $u_1$ is a right Quillen functor,
$u_1X$ is an $\infty$-category
for an $\infty$-bicategory $X$.
The functor $u_1$ induces the functor
$u_1: \bicat\to \cat$
of the underlying $\infty$-categories,
which is a right adjoint to the inclusion functor
$\cat\to \bicat$.

The model structure on 
$\setsc$ is cartesian closed
by \cite[Proposition~3.1.8 and Lemma~4.2.6]{Lurie3}.
(see, also,
\cite[Remark~1.31]{GHL} or
\cite[the paragraph before Lemma~1.24]{GHL2}).
Thus,
we have a function object ${\rm FUN}(X,Y)$
in $\setsc$
for scaled simplicial sets $X$ and $Y$.
When $Y$ is an $\infty$-bicategory,
${\rm FUN}(X,Y)$ is also an $\infty$-bicategory.
We denote by
\[ {\rm Fun}(X,Y) \]
its underlying $\infty$-category
$u_1{\rm FUN}(X,Y)$, and by
\[ \mapsc(X,Y) \]
its underlying $\infty$-groupoid
$u_1{\rm FUN}(X,Y)^{\simeq}$.

The Gray tensor product
of scaled simplicial sets
was introduced in \cite{GHL2}.
This construction determines a functor
\[ (-)\otimes(-):
   \setsc\times\setsc\longrightarrow \setsc, \]
which is a left Quillen bifunctor
by \cite[Theorem~2.14]{GHL2}.
Hence,
it induces a functor
$(-)\otimes (-): \bicat\times\bicat\to\bicat$
of $\infty$-categories,
which preserves colimits separately in each variable.
For example,
the Gray tensor product
$\Delta^1_{\sharp}\otimes \Delta^1_{\sharp}$
is a scaled simplicial set
whose underlying simplicial set
is $\Delta^1\times\Delta^1$ equipped with
the $2$-simplex $\Delta^{\{(0,0),(1,0),(1,1)\}}$
together with degenerate $2$-simplices 
as thin.
It represents an oplax square depicted by
\[ \xymatrix{
  (0,0)\ar[d]\ar[r]&
  (0,1)\ar[d]\\
  (1,0)\ar[r]\ar@{=>}[ur] & (1,1).\\
}\]

The Gray tensor product 
on scaled simplicial sets has adjoints
${\rm FUN}^{\rm lax}(-,-)$ and 
${\rm FUN}^{\rm opl}(-,-)$ which satisfy
\[ \begin{array}{rcl}
   {\rm Hom}_{\setsc}(Z,
   {\rm FUN}^{\rm lax}(X,Y))
   &\cong&
   {\rm Hom}_{\setsc}(Z\otimes X,Y),\\[2mm]
   {\rm Hom}_{\setsc}(Z,
   {\rm FUN}^{\rm opl}(X,Y))
   &\cong&
   {\rm Hom}_{\setsc}(X\otimes Z,Y).\\[2mm]
   \end{array}\]
If $Y$ is an $\infty$-bicategory,
then
${\rm FUN}^{\rm lax}(X,Y)$
and ${\rm FUN}^{\rm opl}(X,Y)$
are also $\infty$-bicategories.
We write
\[ {\rm Fun}^{\rm lax}(X,Y),\quad
   {\rm Fun}^{\rm opl}(X,Y) \]
for the underlying $\infty$-categories
$u_1{\rm FUN}^{\rm lax}(X,Y)$
and $u_1{\rm FUN}^{\rm opl}(X,Y)$,
respectively.

\section{Adjoint morphisms in $(\infty,2)$-categories}
\label{section:adjunctions-infty-bicategory}

In this section
we study adjoint morphisms in $(\infty,2)$-categories.
For an $(\infty,2)$-category $\mathcal{X}$,
we construct $1$-full sub $(\infty,2)$-categories
$\mathcal{X}^L$ and $\mathcal{X}^R$
spanned by left and right adjoint $1$-morphisms,
respectively.
We show that there is a natural dual equivalence
between the underlying $(\infty,1)$-categories
of $\mathcal{X}^L$ and $\mathcal{X}^R$.

First,
we observe a property of adjoint $1$-morphisms
in an $(\infty,2)$-category.
Let $(m^L, m^R, \eta)$ be an adjunction
in an $(\infty,2)$-category $\scriptx$,
where $m^L: x\to x'$ is a left adjoint
$1$-morphism to 
$m^R: x'\to x$,
and $\eta: {\rm id}_{x}\Rightarrow m^Rm^L$ is its unit
$2$-morphism.
For $1$-morphisms
$f: x\to y$ and $g: x'\to y$ of $\scriptx$,
we consider a map of mapping spaces
\[ {\rm Map}_{\scriptx(x',y)}(fm^R,g)\longrightarrow
  {\rm Map}_{\scriptx(x,y)}(f,gm^L) \]
given by the composite
\[ {\rm Map}_{\scriptx(x',y)}(fm^R,g)
   \stackrel{(-)m^L}{\hbox to 10mm{\rightarrowfill}}
   {\rm Map}_{\scriptx(x,y)}(fm^Rm^L,gm^L)
   \stackrel{(f\eta)^*}{\hbox to 10mm{\rightarrowfill}}
   {\rm Map}_{\scriptx(x,y)}(f,gm^L). \]

The following lemma might be well-known
but we give a proof of it for convenience of readers.
   
\begin{lemma}\label{lemma:equivalence-mate-space}
The map
${\rm Map}_{\mathcal{X}(x',y)}(fm^R,g)\to
{\rm Map}_{\mathcal{X}(x,y)}(f,gm^L)$
is an equivalence.
\end{lemma}

\proof
We shall construct a homotopy inverse
of the map.
We consider a map given by the composite
\[ {\rm Map}_{\mathcal{X}(x,y)}(f,gm^L)
    \stackrel{(-)m^R}{\hbox to 10mm{\rightarrowfill}}
    {\rm Map}_{\mathcal{X}(x',y)}(fm^R,gm^Lm^R)
    \stackrel{(g\epsilon)_*}{\hbox to 10mm{\rightarrowfill}}
    {\rm Map}_{\mathcal{X}(x,y)}(fm^R,g),
   \]
where $\epsilon: m^Lm^R\Rightarrow {\rm id}_{x'}$
is a counit of the adjunction. 
By the triangle identity,
we see that
it is the desired homotopy inverse. 
\qed

\bigskip

For an $(\infty,2)$-category $\scriptx$,
we construct 
$1$-full sub $(\infty,2)$-categories $\scriptx^L$ and $\scriptx^R$
of $\scriptx$ spanned by left and right adjoint $1$-morphisms,
respectively.
We give a construction of
$\infty$-bicategories $X^L$ and $X^R$
which represent $\scriptx^L$ and $\scriptx^R$, respectively,
when $\scriptx$
is represented by an $\infty$-bicategory $X$.

Let $X=(\overline{X},T_X)$ be a fibrant scaled simplicial set.
We write $\overline{X}{}^L$
for the subcomplex of $\overline{X}$
spanned by those simplices whose
$1$-dimensional faces are all left adjoint.
We define $X^L$
to be the scaled simplicial set
$(\overline{X}{}^L, T_{X^L})$,
where $T_{X^L}$ is the restriction of
$T_X$ to $\overline{X}{}^L$.
We define a scaled simplicial set
$X^R$ in a similar manner
by using right adjoints
instead of left adjoints.
We show that $X^L$ and $X^R$
are fibrant scaled simplicial sets
and hence they represent
$(\infty,2)$-categories.

\begin{lemma}
The scaled simplicial sets
$X^L$ and $X^R$ are fibrant.
\end{lemma}

\proof
We shall prove the case for $X^L$.
The case for $X^R$ can be proven in a similar manner.
We consider an extension problem of
scaled anodyne maps with respect to $X^L$.
Since $X$ is fibrant,
we can find an extension in $X$.
If a scaled anodyne map is 
of type ({\rm An1}) for $n\ge 3$, ({\rm An2}), or ({\rm An3}),
then the extension lands in $X^L$
since the source of scaled anodyne maps
contains all $1$-simplices of the target.
If a scaled anodyne map is of type ({\rm An1})
for $n=2$,
then the image of $\Delta^{\{0,2\}}$ in $X$
is equivalent to a composite of left adjoints,
and hence the extension also lands in $X^L$.
\qed

\begin{definition}\rm
Let $\scriptx$ be an $(\infty,2)$-category
represented by a fibrant scaled simplicial set $X$.
We denote by $\scriptx^L$ and $\scriptx^R$
the $(\infty,2)$-categories
represented by fibrant
scaled simplicial sets
$X^L$ and $X^R$,
respectively.
\end{definition}

We recall that the inclusion functor
$\cat\to \bicat$ has
the right adjoint $u_1: \bicat\to \cat$.
For an $(\infty,2)$-category $\scriptx$,
$u_1\scriptx$ is the underlying $(\infty,1)$-category
of $\scriptx$
obtained by discarding all non-invertible $2$-morphisms.
Next, we show that
there is a dual equivalence
between $u_1\scriptx^L$ and $u_1\scriptx^R$
by taking adjoint $1$-morphisms.
For this purpose,
we construct a perfect pairing
between $u_1\scriptx^L$
and $u_1\scriptx^R$.

First,
we construct a pairing
of $u_1\scriptx^L$ and $u_1\scriptx^R$.
For this purpose,
we recall the oplax arrow $\infty$-category
\[ \oplaxarrow{\scriptx} \]
which is represented by
the $\infty$-category
${\rm Fun}^{\rm opl}(\Delta^1_{\sharp},X)$
(cf.~\cite[Definition~7.5]{HHLN2}).
An object of $\oplaxarrow{\scriptx}$
is a $1$-morphism of $\scriptx$.
A morphism of $\oplaxarrow{\scriptx}$
from $f: x\to y$ to $f': x'\to y'$
is an oplax square
depicted by
\[ \xymatrix{
  x\ar[d]_{f}\ar[r]&
  x'\ar[d]^{f'}\\
  y\ar[r]\ar@{=>}[ur] & y'.\\
}\]

The inclusion map $\partial \Delta^{1}_{\sharp}
\to \Delta^{1}_{\sharp}$
induces a map
${\rm FUN}^{\rm opl}(\Delta^1_{\sharp},X)\to
X\times X$,
which is a fibration between fibrant objects
in $\setsc$.
Since $u_1: \setsc\to \setjoy$ is a right Quillen
functor,
we obtain a map
$\oplaxarrow{X}\to u_1X\times u_1X$,
which is a fibration of fibrant objects
in $\setjoy$.
We write
\[ (s,t): \oplaxarrow{\scriptx}\longrightarrow
   u_1\scriptx\times u_1\scriptx \]
for the corresponding map of $(\infty,1)$-categories.
          
\begin{remark}\rm
The map
$(s,t):\oplaxarrow{\scriptx}\to
u_1\scriptx\times u_1\scriptx$
is an orthofibration
by \cite[Proposition~7.9]{HHLN2}.
\end{remark}


Now,
we give a sufficient condition for $1$-morphisms 
in $\oplaxarrow{\scriptx}$
to be $(s,t)$-coCartesian.

\begin{lemma}\label{lemma:oplax-arrow-L-cocart}
Let $\sigma$ be a morphism
in $\oplaxarrow{\scriptx}$ depicted by 
\[ \xymatrix{
  x\ar[d]_{f}\ar[r]^{g^L} &
  x'\ar[d]^{f'}\\
  y\ar[r]_h\ar@{=>}[ur]^u & y'\\}\]
where $g^L$ is a left adjoint $1$-morphism.
Then $\sigma$ is an $(s,t)$-coCartesian morphism
if $f'$ is equivalent to $hfg^R$,
and 
the $2$-morphism $u$
is equivalent to
the composite
$hf\stackrel{hf\eta}{\Longrightarrow}
hfg^Rg^L\simeq f'g^L$,
where $(g^L,g^R,\eta)$
is an adjunction in $\scriptx$
that is an extension of $g^L$.
\end{lemma}

\proof
Let $f'': x''\to y''$ be an object
of $\oplaxarrow{\scriptx}$.
We would like to show that
the following commutative diagram
\[ \begin{array}{ccc}
    {\rm Map}_{\oplaxarrow{\scriptx}}(f',f'')
    & \longrightarrow &
    {\rm Map}_{\oplaxarrow{\scriptx}}(f,f'')  \\[1mm]
    \bigg\downarrow & & \bigg\downarrow \\[3mm]
    {\rm Map}_{u_1\scriptx\times u_1\scriptx}
    ((x',y'),(x'',y''))
    & \longrightarrow &
    {\rm Map}_{u_1\scriptx\times u_1\scriptx}
    ((x,y),(x'',y''))\\
\end{array}\]
is a pullback in
the $(\infty,1)$-category of spaces.
We take a $1$-morphism
$(k,l): (x',y')\to (x'',y'')$
in $u_1\scriptx\times u_1\scriptx$.
The induced map on fibers
\[ \map{\oplaxarrow{\scriptx}}
   (f',f'')_{(k,l)}
   \to
   \map{\oplaxarrow{\scriptx}}
   (f,f'')_{(kg^L,lh)} \]
is identified with the map
\[ \map{\scriptx(x',y'')}
   (lf',f''k)
   \longrightarrow 
   \map{\scriptx(x,y'')}
   (lhf,f''kg^L) \]
given by the composition with $\sigma$, 
which is an equivalence by Lemma~\ref{lemma:equivalence-mate-space}.
\qed

\bigskip

We define maps
\[ \begin{array}{rrcl}
  (s,t)_{sL}:& \oplaxarrowsl{\scriptx}
  &\longrightarrow&
  u_1\scriptx^L\times u_1\scriptx,\\[2mm]
  (s,t)_{L,R}:& \oplaxarrowLR{\scriptx}&
  \longrightarrow&
  u_1\scriptx^L\times u_1\scriptx^R \\
  \end{array}\]
by the following pullback diagrams
\[ \begin{array}{ccccc}
    \oplaxarrowLR{\scriptx}
     & \longrightarrow &
     \oplaxarrowsl{\scriptx}  
    & \longrightarrow &
     \oplaxarrow{\scriptx}\\     
    \mbox{$\scriptstyle (s,t)_{L,R}$}
    \bigg\downarrow
    \phantom{\mbox{$\scriptstyle (s,t)_{L,R}$}}& &
    \mbox{$\scriptstyle (s,t)_{sL}$}
    \bigg\downarrow
    \phantom{\mbox{$\scriptstyle (s,t)_{sL}$}}& &
    \phantom{\mbox{$\scriptstyle (s,t)$}}
    \bigg\downarrow
    \mbox{$\scriptstyle (s,t)$}\\
    u_1\scriptx^L\times u_1\scriptx^R
    & \longrightarrow &
    u_1\scriptx^L\times u_1\scriptx
    & \longrightarrow &
    u_1X\times u_1X.
   \end{array}\]
By the criterion of
$(s,t)$-coCartesian morphisms
in Lemma~\ref{lemma:oplax-arrow-L-cocart},
we obtain the following corollary.

\begin{corollary}\label{cor:Ar-sL-coCart}
The maps
$(s,t)_{sL}:
\oplaxarrowsl{\scriptx}
\to u_1\scriptx^L\times u_1\scriptx$
and
$(s,t)_{L,R}:
\oplaxarrowLR{\scriptx}
\to u_1\scriptx^L\times u_1\scriptx^R$
are coCartesian fibrations.
\end{corollary}

Next,
we construct the desired pairing
from the map $(s,t)_{L,R}$.
We define
$\overlineoplaxarrow{\scriptx}$
to be the full subcategory
of $\oplaxarrowLR{\scriptx}$
spanned by those objects that are
right adjoint $1$-morphisms in $\scriptx$.
By restricting the map $(s,t)_{L,R}$ to
$\overlineoplaxarrow{\scriptx}$,
we obtain a map
\[ \overline{(s,t)}_{L,R}:
\overlineoplaxarrow{\scriptx}
\longrightarrow u_1\scriptx^L\times u_1\scriptx^R. \]

\begin{lemma}\label{lemma:overline-arrow-coCart}
  The map $\overline{(s,t)}_{L,R}:
  \overlineoplaxarrow{\scriptx}
  \to u_1\scriptx^L\times u_1\scriptx^R$
is a coCartesian fibration.
\end{lemma}

\proof
Let $f: x\to y$ be an object of
$\overlineoplaxarrow{\scriptx}$
and let
$(g^L,h^R): (x,y)\to (x',y')$ be a morphism
of $u_1\scriptx^L\times u_1\scriptx^R$.
We take a $(s,t)_{L,R}$-coCartesian morphism
$f\to f'$ in
$\oplaxarrowLR{\scriptx}$
covering $(g^L,h^R)$.
It suffices to show that
$f'$ belongs to $\overlineoplaxarrow{\scriptx}$.
This follows from the fact that $h^Rfg^R$ is
a right adjoint $1$-morphism.
\qed

\bigskip

We define ${\rm Pair}(\scriptx)$
to be the wide subcategory of
$\overlineoplaxarrow{\scriptx}$
spanned by $\overline{(s,t)}_{L,R}$-coCartesian morphisms.
We denote by
\[ (s,t)^{\rm pair}: {\rm Pair}(\scriptx)
   \longrightarrow u_1\scriptx^L\times u_1\scriptx^R \]
the restriction of $\overline{(s,t)}_{L,R}$
to ${\rm Pair}(\scriptx)$.
By Lemma~\ref{lemma:overline-arrow-coCart},
the map $(s,t)^{\rm pair}$
is a left fibration,
and hence it defines a pairing of $(\infty,1)$-categories
between $u_1\scriptx^L$ and $u_1\scriptx^R$.

\begin{theorem}
\label{theorem:L-R-perfect-pairing}
The pairing 
$(s,t)^{\rm pair}:{\rm Pair}(\scriptx)
\to u_1\scriptx^L\times u_1\scriptx^R$
is perfect.
\end{theorem}

\proof
By \cite[Corollary~5.2.1.22]{Lurie2},
it suffices to show that
the pairing $(s,t)^{\rm pair}$
is both left and right representable
and that an object of ${\rm Pair}(\scriptx)$
is left universal if and only if it is right universal.

First,
we show that $(s,t)^{\rm pair}$ is left representable.
Let $x$ be an object of $u_1\scriptx^L$.
We set ${\rm S}(x)=\{x\}\times_{u_1\scriptx^L}
{\rm Pair}(\scriptx)$
for simplicity.
We take the identity morphism
${\rm id}_x: x\to x$,
which we regard as an object of ${\rm S}(x)$.
We have to show that ${\rm id}_x$
is an initial object.
For any object $f: x\to y$
in ${\rm S}(x)$,
we can identify
$\map{{\rm S}(x)}({\rm id}_x,f)$
with $\map{u_1\scriptx^R}(x,y)_{/f}$,
which is contractible.

Next,
we shall show that $(s,t)^{\rm pair}$ is right representable.
Let $y$ be an object of $u_1\scriptx^R$.
We set ${\rm T}(y)={\rm Pair}(\scriptx)
\times_{u_1\scriptx^R}\{y\}$
for simplicity.
We take the identity morphism
${\rm id}_y: y\to y$,
which we regard as an object of ${\rm T}(y)$.
We have to show that ${\rm id}_y$
is an initial object.
For any object $f: x\to y$
in ${\rm T}(y)$,
we can identify
$\map{{\rm T}(y)}({\rm id}_y,f)$
with the space of units $(g,f,u)$
that are extensions of $f$,
which is contractible
by \cite[Theorem~4.4.18]{Riehl-Verity}.

Finally,
we shall show that
an object of ${\rm Pair}(\scriptx)$
is left universal if and only if
it is right universal.
By the above argument,
we see that an object $f: x\to y$ of
${\rm Pair}(\scriptx)$
is left universal if and only if
$f$ is equivalent to ${\rm id}_x$
in ${\rm S}(x)$.
This means that $f$ is an equivalence in $\scriptx$.
In a similar manner
we see that $f$ is right universal if and only
if $f$ is an equivalence in $\scriptx$.
This completes the proof.
\qed

\begin{corollary}
\label{cor:dual-eq-u1-Z-L-R}
There is a dual equivalence of $(\infty,1)$-categories
between $u_1\scriptx^L$ and $u_1\scriptx^R$.
\end{corollary}

\begin{remark}\rm
The dual equivalence in Corollary~\ref{cor:dual-eq-u1-Z-L-R}
is the identity on objects 
and
assigns to a left adjoint $1$-morphism
its right adjoint $1$-morphism.
\end{remark}

Finally,
we show that the dual equivalence
in Corollary~\ref{cor:dual-eq-u1-Z-L-R}
is natural with respect to $\scriptx$.
We let $u_1(-)^L$ and $u_1(-)^R$
be functors from $\bicat$ to $\cat$
which assigns $u_1\scripty^L$ and
$u_1\scripty^R$ to
an $(\infty,2)$-category $\scripty$,
respectively.

\begin{corollary}
\label{cor:functor-perfect-pairing-LR}
There is a natural dual equivalence of functors
between $u_1(-)^L$ and $u_1(-)^R$.
\end{corollary}

\proof
Let $f: \scripty\to \scriptz$
be a functor of $(\infty,2)$-categories.
This induces a morphism of pairings
\begin{equation}\label{eq:morphism-pairing}
  \begin{array}{ccc}
  {\rm Pair}(\scripty)
  &\stackrel{{\rm Pair}(f)}{\hbox to 20mm{\rightarrowfill}}&
  {\rm Pair}(\scriptz)\\
  \bigg\downarrow & & \bigg\downarrow\\
  u_1\scripty^L\times u_1\scripty^R 
  & \stackrel{u_1f^L\times u_1f^R}
  {\hbox to 20mm{\rightarrowfill}}&
  u_1\scriptz^L\times u_1\scriptz^R.\\   
  \end{array}
  \end{equation}
By the characterization
of right universal objects
in ${\rm Pair}(\scriptx)$
in the proof of Theorem~\ref{theorem:L-R-perfect-pairing},
we see that (\ref{eq:morphism-pairing})
is right representable
in the sense of \cite[Variant~5.2.1.16]{Lurie2}.
By this observation,
we obtain a functor
\[ \bicat\longrightarrow {\rm CPair}^{\rm perf} \]
of $\infty$-categories
by assigning to an $(\infty,2)$-category $\scripty$
the perfect pairing ${\rm Pair}(\scripty)\to
u_1\scripty^L\times u_1\scripty^R$,
where ${\rm CPair}^{\rm perf}$
is the $\infty$-category of perfect pairings
and right representable morphisms
(see \cite[Remark~5.2.1.20]{Lurie2}
for the definition of ${\rm CPair}^{\rm perf}$).
This implies that
the dual equivalence
in Corollary~\ref{cor:dual-eq-u1-Z-L-R}
is natural with respect to
$\scriptx$.
\qed

\begin{remark}\rm
We can prove an existence of a natural dual equivalence between
$u_1\scriptx^L$ and $u_1\scriptx^R$
by using \cite[Remark~4.11]{Haugseng4}
as well.
Let $\mathfrak{adj}$ be the free $2$-category
containing an adjunction,
and let $\mathfrak{adj}^{\bullet}$
be the cosimplicial $(\infty,2)$-category
whose $n$th object
is the $(\infty,2)$-category of $n$ composable
adjunctions
$\mathfrak{adj}\coprod_{[0]}\cdots
\coprod_{[0]}\mathfrak{adj}$.
By \cite[Theorem~4.4 and Remark~4.11]{Haugseng4},
there are equivalences
of simplicial spaces
\[ {\rm Map}_{\bicat}([\bullet]^{\rm op},\scriptx)^{\rm ladj}
   \stackrel{\simeq}{\longleftarrow}
   {\rm Map}_{\bicat}
   (\mathfrak{adj}^{\bullet},\scriptx)
   \stackrel{\simeq}{\longrightarrow}
            {\rm Map}_{\bicat}([\bullet],\scriptx)^{\rm radj}.\]
This implies
that there is a natural dual equivalence
of complete Segal spaces
associated to $u_1\scriptx^L$
and $u_1\scriptx^R$.
\end{remark}

\section{Map monoidales in monoidal $\infty$-bicategories}
\label{section:map-monoidale}

In this section
we generalize map monoidales in monoidal bicategories
in the setting of higher category theory.
We introduce notions of left adjoint $\scripto$-monoidales
and right adjoint $\scripto$-comonoidales
in an $\scripto$-monoidal $(\infty,2)$-category.
We show that there is a dual equivalence
between the $(\infty,1)$-category of left adjoint
$\scripto$-monoidales and
the $(\infty,1)$-category
of right adjoint $\scripto$-comonoidales 
by taking adjoint morphisms.



Since $\bicat$ has finite products,
we can consider $\scripto$-monoid objects
in $\bicat$.
We call an $\scripto$-monoid object
in $\bicat$
an $\scripto$-monoidal $(\infty,2)$-category.  
In the following of this section
we let $\scriptz$
be an $\scripto$-monoidal $(\infty,2)$-category.

Since the functor $u_1: \bicat\to\cat$
is right adjoint to the inclusion functor,
it preserves small limits.
In particular,
it preserves finite products
and hence it induces a functor
$u_1: {\rm Mon}_{\scripto}(\bicat)\to
{\rm Mon}_{\scripto}(\cat)$.
From this observation,
we see that $u_1\scriptz$ is
an $\scripto$-monoidal $(\infty,1)$-category.

Let $\scriptw$ be 
an $\scripto$-monoidal $(\infty,1)$-category.
We identify $\scriptw$
with a coCartesian fibration
$\scriptw^{\otimes}\to\scripto^{\otimes}$
of $\infty$-operads.
An $\scripto$-algebra object $A$
in $\scriptw$ is a map
of $\infty$-operads $\scripto^{\otimes}
\to \scriptw^{\otimes}$ over $\scripto^{\otimes}$.
We define 
${\rm Alg}_{\scripto}(\scriptw)$
to be the full subcategory of
${\rm Fun}_{\scripto^{\otimes}}(\scripto^{\otimes},
\scriptw^{\otimes})$
spanned by $\scripto$-algebra objects,
and call it 
the $\infty$-category of $\scripto$-algebra
objects in $\scriptw$.

The opposite $(\infty,1)$-category $\scriptw^{\rm op}$
inherits an $\scripto$-monoidal structure from $\scriptw$.
An $\scripto$-coalgebra object in $\scriptw$ is
an $\scripto$-algebra object in $\scriptw^{\rm op}$.
We define
${\rm coAlg}_{\scripto}(\scriptw)$
to be
${\rm Alg}_{\scripto}(\scriptw^{\rm op})^{\rm op}$,
and call it 
the $\infty$-category of $\scripto$-coalgebra
objects in $\scriptw$.

In ordinary category theory
a monoidale is a pseudomonoid
in a monoidal bicategory. 
We define a monoidale in the setting
of higher category theory.

\begin{definition}\rm
Let $\scriptz$ be an $\scripto$-monoidal
$(\infty,2)$-category.
An $\scripto$-monoidale in $\mathcal{Z}$
is an $\scripto$-algebra object in $u_1\scriptz$.
An $\scripto$-comonoidale in $\scriptz$ 
is an $\scripto$-coalgebra object in $u_1\scriptz$.
\end{definition}

Next,
we would like to introduce a notion
of a map monoidale in the setting of higher category theory.
For this purpose,
we show that
the $(\infty,2)$-categories
$\scriptz^L$ and $\scriptz^R$
inherit $\scripto$-monoidal structure from $\scriptz$.

\begin{lemma}\label{lemma:L-R-monoidal-structure-support}
If $\scriptz$ is an $\scripto$-monoidal
$(\infty,2)$-category,
then $\scriptz^L$ and $\scriptz^R$
are also $\scripto$-monoidal $(\infty,2)$-categories.
\end{lemma}

\proof
We shall show that $\scriptz^L$
is an $\scripto$-monoidal $(\infty,2)$-category.
The case of $\scriptz^R$ is similar.
The construction $\scriptx\mapsto \scriptx^L$
determines a functor $(-)^L: \bicat\to\bicat$.
This functor preserves finite products.
Thus,
it induces a functor
$(-)^L: {\rm Mon}_{\scripto}(\bicat)\to {\rm Mon}_{\scripto}(\bicat)$.
\qed

\bigskip

By Lemma~\ref{lemma:L-R-monoidal-structure-support},
we have
$\scripto$-monoidal $(\infty,1)$-categories
$u_1\scriptz^L$ and $u_1\scriptz^R$.
In ordinary category theory
a monoidale is said to be a map monoidale
if the product and unit morphisms
are left adjoint $1$-morphisms.
We define a notion of a map monoidale
in an $\scripto$-monoidal $(\infty,2)$-category.  

\begin{definition}\rm
Let $\scriptz$ be an $\scripto$-monoidal
$(\infty,2)$-category.
We define a map $\scripto$-monoidale in $\scriptz$
to be an $\scripto$-monoidale in $\scriptz^L$.
We also say that
a map $\scripto$-monoidale
is a left adjoint $\scripto$-monoidale.
We define a right adjoint $\scripto$-comonoidale in $\scriptz$
to be an $\scripto$-comonoidale in $\scriptz^R$.
\end{definition}


Next,
we study a duality of
left adjoint $\scripto$-monoidales and
right adjoint $\scripto$-comonoidales
under the dual equivalence
in Corollary~\ref{cor:dual-eq-u1-Z-L-R}.
For this purpose,
we promote the dual equivalence
in Corollary~\ref{cor:dual-eq-u1-Z-L-R}
to a dual equivalence of monoidal
$(\infty,1)$-categories.

\begin{lemma}
\label{lemma:monoidal-dual-equivalence-LR}
When $\scriptz$ is an $\scripto$-monoidal
$(\infty,2)$-category,
there is a natural dual equivalence
of $\scripto$-monoidal $(\infty,1)$-categories
between $u_1\scriptz^L$ and $u_1\scriptz^R$.
\end{lemma}

\proof
By the proof of Corollary~\ref{cor:functor-perfect-pairing-LR},
there is a functor 
$\bicat\to {\rm CPair}^{\rm pair}$
which assigns to an $(\infty,2)$-category
$\scriptx$ the perfect pairing
${\rm Pair}(\scriptx)\to u_1\scriptx^L\times
\scriptx^R$.
We can verify that
this functor preserves finite products.
Hence we obtain a functor
${\rm Mon}_{\scripto}(\bicat)\to
{\rm Mon}_{\scripto}({\rm Pair}^{\rm perf})$.
This implies that 
the dual equivalence in Corollary~\ref{cor:dual-eq-u1-Z-L-R}
can be promoted to a natural dual equivalence
of $\scripto$-monoidal $(\infty,1)$-categories.
\qed

\bigskip

By Lemma~\ref{lemma:monoidal-dual-equivalence-LR},
we obtain a duality
between left adjoint $\scripto$-monoidales 
and right adjoint $\scripto$-comonoidales.

\begin{corollary}\label{cor:alg-coalg-correspondence}
There is a natural dual equivalence between
the $\infty$-category of left adjoint $\scripto$-monoidales
and the $\infty$-category of right adjoint
$\scripto$-comonoidales
\[ {\rm Alg}_{\scripto}(u_1\scriptz^L)\simeq 
   {\rm coAlg}_{\scripto}(u_1\scriptz^R)^{\rm op}. \]
\end{corollary}

\begin{definition}\rm
For a left adjoint $\scripto$-monoidale
$A\in {\rm Alg}_{\scripto}(u_1\scriptz^L)$,
we let $A^*$
be the corresponding object
in ${\rm coAlg}_{\scripto}(u_1\scriptz^R)$
under the equivalence
in Corollary~\ref{cor:alg-coalg-correspondence}.
We call $A^*$ the right adjoint $\scripto$-comonoidale
associated to $A$.
Similarly,
for a right adjoint $\scripto$-comonoidale
$C\in {\rm coAlg}_{\scripto}(u_1\scriptz^R)$,
we denote by $C^*$
the corresponding object
in ${\rm Alg}_{\scripto}(u_1\scriptz^L)$
under the equivalence
in Corollary~\ref{cor:alg-coalg-correspondence}.
We call $C^*$ the left adjoint $\scripto$-monoidale
associated to $C$.
\end{definition}

\section{Map monoidales and Duoidal $\infty$-categories}
\label{section:map-duoidal}

Let $A$ be a map $\scripto$-monoidale
in an $\scripto$-monoidal
$(\infty,2)$-category $\scriptz$.
The goal of this section 
is to show that the endomorphism mapping
$\infty$-category $\scriptz(A,A)$
admits the structure of
a coCartesian $(\Delta^{\rm op},\scripto)$-duoidal
$\infty$-category.

First,
we construct a
simplicial $\scripto$-monoidal $(\infty,1)$-category
$\funopl([\bullet],\scriptz)^{\otimes}$.
For an $(\infty,2)$-category $\mathcal{K}$,
the functor $\funopl(\mathcal{K},-):
\bicat\to\cat$ preserves finite products.
Hence,
it induces a functor
\[ \funopl(\mathcal{K},-)^{\otimes}:
{\rm Mon}_{\scripto}(\bicat)\to {\rm Mon}_{\scripto}(\cat). \]
Thus,
$\funopl
(\mathcal{K},\scriptz)^{\otimes}$ is an
$\scripto$-monoidal $(\infty,1)$-category
for an $\scripto$-monoidal $(\infty,2)$-category
$\scriptz$.
We notice that
it is functorial with respect to $\mathcal{K}$.
We define a simplicial $\scripto$-monoidal
$(\infty,1)$-category
\[ \funopl([\bullet],\scriptz)^{\otimes}:
   \Delta^{\rm op}\longrightarrow
   {\rm Mon}_{\scripto}(\cat)   \]
which is given by
$[n]\mapsto \funopl([n],\scriptz)^{\otimes}$.

Notice that
the $0$th $\scripto$-monoidal $(\infty,1)$-category
$\funopl([0],\scriptz)^{\otimes}$
is equivalent to $u_1\scriptz^{\otimes}$.
The inclusion map
$i: \{[0]\}\hookrightarrow \Delta^{\rm op}$
induces an adjunction
\[ i^*: {\rm Fun}(\Delta^{\rm op},{\rm Mon}_{\scripto}(\cat))
        \rightleftarrows
            {\rm Mon}_{\scripto}(\cat): i_*,\]
where $i^*$ is the restriction
and $i_*$ is its right adjoint
given by the right Kan extension along $i$.
Hence,
the equivalence
$i^*\funopl([\bullet],\scriptz)^{\otimes}\simeq
u_1\scriptz^{\otimes}$
gives rise to a map
\[ \funopl([\bullet],\scriptz)^{\otimes}
   \longrightarrow i_*u_1\scriptz^{\otimes} \]
in 
${\rm Fun}(\Delta^{\rm op},{\rm Mon}_{\scripto}(\cat))$.
By using the composite
${\rm Mon}_{\scripto}(\cat)\simeq
\mathsf{Mon}_{\scripto}(\cat)\to
\mathsf{Mon}_{\scripto}^{\rm lax}(\cat)
\to \op{\scripto^{\otimes}}$,
we regard it as a map
in ${\rm Fun}(\Delta^{\rm op},\op{\scripto^{\otimes}})$.

Let $A\in {\rm Alg}_{\scripto}(u_1\scriptz)$
be an $\scripto$-monoidale in $\scriptz$.
We regard $A$ as a morphism
of $\infty$-operads 
$A: \scripto^{\otimes}\to
      u_1\scriptz^{\otimes}$
over $\scripto^{\otimes}$.
By the right Kan extension along $i$,
we obtain a morphism
\[ i_*A: i_*\scripto^{\otimes}
         \longrightarrow i_*u_1\scriptz^{\otimes} \]
in ${\rm Fun}(\Delta^{\rm op},
{\rm Op}_{\infty/\scripto^{\otimes}})$.

\begin{definition}\rm
We define $\scriptz(A,A)^{\otimes}_{\bullet}$
by the following pullback diagram
in ${\rm Fun}(\Delta^{\rm op},
{\rm Op}_{\infty/\scripto^{\otimes}})$
\begin{equation}\label{eq:defining-diagram-ZAA}
  \begin{array}{ccc}
    \scriptz(A,A)^{\otimes}_{\bullet}&
     \longrightarrow &
     \funopl([\bullet],\scriptz)^{\otimes}\\
     \bigg\downarrow &&\bigg\downarrow \\
     i_*\scripto^{\otimes}&
     \stackrel{i_*A}{\longrightarrow}&
     i_*u_1\scriptz^{\otimes}.\\
\end{array}
\end{equation}
\end{definition}

We would like to show that
$\scriptz(A,A)^{\otimes}_{\bullet}$
is a $\Delta^{\rm op}$-monoid object
in $\mathsf{Mon}_{\scripto}^{\rm lax}(\cat)$,
which is a full subcategory of $\op{\scripto^{\otimes}}$,
when $A$ is a map $\scripto$-monoidale.
First,
we show that $\scriptz(A,A)^{\otimes}_{\bullet}$
is a $\Delta^{\rm op}$-monoid object
of $\op{\scripto^{\otimes}}$,
and hence it corresponds to 
a virtual coCartesian $(\Delta^{\rm op},\scripto)$-duoidal
$\infty$-category
by Remark~\ref{remark:cor-virtual-duoidal-category}.

\begin{proposition}\label{prop:Z-A-A-virtual-duoidal}
The functor $\scriptz(A,A)^{\otimes}_{\bullet}$
is a $\Delta^{\rm op}$-monoid object
of $\op{\scripto^{\otimes}}$.
\end{proposition}

\proof
For a simplicial object $F$
in an $(\infty,1)$-category
which admits finite limits,
we say that $F$ is a category object if 
the inert maps $[n]\to [i]$ in $\Delta^{\rm op}$
for $i=0,1$
induce an equivalence
\[ F([n])\to F(\{0,1\})\times_{F(\{1\})}\cdots
  \times_{F(\{n-1\})}F(\{n-1,n\}) \]
for each $n\ge 0$.
In order to prove the proposition,
it suffices to show that 
$\scriptz(A,A)^{\otimes}_{\bullet}$
is a category object
and that $\scriptz(A,A)^{\otimes}_0$
is a final object in $\op{\scripto^{\otimes}}$.

First,
we shall show that 
$\scriptz(A,A)^{\otimes}_{\bullet}$
is a category object
in $\op{\scripto^{\otimes}}$.
Since category objects
are closed under limits,
it suffices to show that
$i_*\scripto^{\otimes}$,
$i_*u_1\scriptz^{\otimes}$,
and
$\funopl([\bullet],\scriptz)^{\otimes}$
are category objects.
For an $\infty$-operad
$\scriptp^{\otimes}$ over $\scripto^{\otimes}$,
the simplicial object $i_*\scriptp^{\otimes}$
is given by
\[ [n]\mapsto
   \overbrace{\scriptp^{\otimes}\times_{\scripto^{\otimes}}
     \cdots\times_{\scripto^{\otimes}}\scriptp^{\otimes}}^{n+1}.\]
By this description,
we see that $i_*\scripto^{\otimes}$ and
$i_*u_1\scriptz^{\otimes}$ are category objects.
By using an equivalence
$\{0,1\}\coprod_{\{1\}}
   \cdots \coprod_{\{n-1\}}\{n-1,n\}
   \to [n]$
of $(\infty,1)$-categories
and the fact that the Gray tensor product
preserves colimits in each variable,
we obtain an equivalence

\[  (\{0,1\}{\otimes}\mathcal{Y})
    \coprod_{\{1\}{\otimes}\mathcal{Y}}
    \cdots
    \coprod_{\{n-1\}{\otimes}\mathcal{Y}}
    (\{n-1,n\}{\otimes}\mathcal{Y})
    \stackrel{\simeq}{\longrightarrow}
    [n]{\otimes}\mathcal{Y}\]
of $(\infty,2)$-categories   
for any $(\infty,1)$-category $\mathcal{Y}$.
This equivalence
implies that
the simplicial object
$\map{\cat}(\mathcal{Y},
\funopl([\bullet],\scriptz)^{\otimes})$
is a category object
in ${\rm Mon}_{\scripto}(\mathcal{S})$.
By the Yoneda lemma
and the fact that
${\rm Mon}_{\scripto}(\cat)\to\op{\scripto^{\otimes}}$
preserves limits,
we see that $\funopl([\bullet],\scriptz)^{\otimes}$
is a category object of
$\op{\scripto^{\otimes}}$.

Next,
we shall show that
$\scriptz(A,A,)^{\otimes}_0$
is a final object of $\op{\scripto^{\otimes}}$.
Since the right vertical arrow in
(\ref{eq:defining-diagram-ZAA})
induces an equivalence,
so does the left vertical arrow.
The claim follows from
the equivalence
$i_*\scripto^{\otimes}([0])\simeq
\scripto^{\otimes}$.
\qed

\bigskip

Finally,
when $A$ is a map
$\scripto$-monoidale in $\scriptz$,
we show that
$\scriptz(A,A)^{\otimes}_{\bullet}$
is a $\Delta^{\rm op}$-monoid object
of $\mathsf{Mon}_{\scripto}^{\rm lax}(\cat)$,
which corresponds to a 
coCartesian $(\Delta^{\rm op},\scripto)$-duoidal
$\infty$-category
by Remark~\ref{remark:cor-duoidal-category}.



\begin{theorem}\label{thm:Z-A-A-duoidal-structure}
We assume that $A$ is a map $\scripto$-monoidale
in an $\scripto$-monoidal $(\infty,2)$-category $\scriptz$.   
Then,
the functor $\scriptz(A,A)^{\otimes}_{\bullet}$
is an object
of ${\rm Mon}_{\Delta^{\rm op}}
(\mathsf{Mon}_{\scripto}^{\rm lax}(\cat))$,
and hence
it corresponds to a coCartesian
$(\Delta^{\rm op},\scripto)$-duoidal $\infty$-category.
\end{theorem}

\proof
By Proposition~\ref{prop:Z-A-A-virtual-duoidal},
it suffices to show that
$\scriptz(A,A)^{\otimes}_1$ is
an $\scripto$-monoidal $\infty$-category.
We write
$\oplaxarrow{\scriptz}^{\otimes}$ 
for $\funopl([1],\scriptz)^{\otimes}$.
Since $A$ is a map $\scripto$-monoidale,
the map $A: \scripto^{\otimes}\to u_1\scriptz^{\otimes}$
factors through $(u_1\scriptz^{L})^{\otimes}$.
We consider the following pullback diagram
\[ \begin{array}{ccc}
     \oplaxarrowsl{\scriptz}
      & \longrightarrow &
      \oplaxarrow{\scriptz}^{\otimes} \\
     \mbox{$\scriptstyle (s,t)_{sL}^{\otimes}$}
     \bigg\downarrow
     \phantom{\mbox{$\scriptstyle (s,t)_{L,L}^{\otimes}$}}
     & &
     \phantom{\mbox{$\scriptstyle (s,t)^{\otimes}$}}
     \bigg\downarrow
     \mbox{$\scriptstyle (s,t)^{\otimes}$}\\
     (u_1\scriptz^{L})^{\otimes}\times_{\scripto^{\otimes}}
     u_1\scriptz^{\otimes}
     & \longrightarrow &
     u_1\scriptz^{\otimes}\times_{\scripto^{\otimes}}
     u_1\scriptz^{\otimes} \\
   \end{array}\]
in ${\rm Op}_{\infty/\scripto^{\otimes}}$.
Since $\scriptz(A,A)^{\otimes}_1\to\scripto^{\otimes}$
is obtained by taking pullback of
$(s,t)^{\otimes}_{sL}$ along the map
$(A,A): \scripto^{\otimes}\to
(u_1\scriptz^{L})^{\otimes}\times_{\scripto^{\otimes}}
u_1\scriptz^{\otimes}$,
it suffices to show that
$(s,t)^{\otimes}_{sL}$ is a coCartesian fibration. 

We have a commutative diagram
\[ \xymatrix{
  \oplaxarrowsl{\scriptz}
  \ar[rr]^-{(s,t)_{sL}^{\otimes}}\ar[dr]_p& &
   (u_1\scriptz^{L})^{\otimes}\times_{\scripto^{\otimes}}
   u_1\scriptz^{\otimes}\ar[dl]^q\\
  & \ \scripto^{\otimes},& \\ 
}\]
where $p$ and $q$ are coCartesian fibrations.
First,
we shall show that
$(s,t)_{sL}^{\otimes}$ is a locally coCartesian fibration.
For this,
by the dual of \cite[Proposition~2.4.2.11]{Lurie1},
it suffices to show that the following two conditions
hold:
(1)
The map $(s,t)_{sL}^{\otimes}$ carries
$p$-coCartesian morphisms
to $q$-coCartesian morphisms.
(2)
For any $a\in \scripto^{\otimes}$,
the induced map 
$\oplaxarrowsl{\scriptz}_{a}
\to (u_1\scriptz^{L})^{\otimes}_a
\times u_1\scriptz^{\otimes}_a$
on fibers
is a coCartesian fibration.

For (1), 
we let $f: x\to x'$ be an object
of $\oplaxarrowsl{\scriptz}$
over $a\in \scripto^{\otimes}$,
and
let $\phi: a\to b$ be a morphism in $\scripto^{\otimes}$.
We may assume that
$\phi:a \to b$ is active,
where $a\simeq \{a_i\}_{i\in |I|}$,
and $a_i,b\in\scripto$.
We can write
$f\simeq \{f_i\}_{i\in |I|}$,
where $f_i: x_i\to x_i'$.
A $p$-coCartesian morphism
$f\to \phi_*f$ is identified with
$\{f_i\}_{i\in |I|}\to \otimes_{\phi}f_i$.
Then 
$(s,t)_{sL}^{\otimes}(f\to \phi_*f)$
is identified with
$\{(x_i,x_i')\}_{i\in |I|}\to
(\otimes_{\phi}x_i,\otimes_{\phi}x_i')$.
Thus,
$(s,t)_{sL}^{\otimes}$ carries $p$-coCartesian morphisms
to $q$-coCartesian morphisms.
For (2),
we may assume that $a\in\scripto$.
Then the map
$\oplaxarrowsl{\scriptz}_{a}
\to (u_1\scriptz^{L})^{\otimes}_a
\times u_1\scriptz^{\otimes}_a$
is a coCartesian fibration
by Corollary~\ref{cor:Ar-sL-coCart}.

Finally,
we would like to show that
$(s,t)_{sL}^{\otimes}$
is in fact a coCartesian fibration.
By the dual of \cite[Proposition~2.4.2.8]{Lurie1}
together with the characterization of
locally $(s,t)_{sL}^{\otimes}$-coCartesian morphisms
in the dual of \cite[Proposition~2.4.2.11]{Lurie1},
it suffices to show that the following condition holds:
(3)
For any morphism
$\phi: a\to b$ in $\scripto^{\otimes}$,
the induced functor
$\phi_*:
\oplaxarrowsl{\scriptz}_a
\to
\oplaxarrowsl{\scriptz}_b$
carries $(s,t)_{sL,a}^{\otimes}$-coCartesian morphisms
to $(s,t)_{sL,b}^{\otimes}$-coCartesian morphisms.

For (3),
we may assume that
$\phi:a \to b$ is active,
where $a\simeq \{a_i\}_{i\in |I|}$,
and $a_i,b\in\scripto$.
We let
$\{f_i\}_{i\in |I|}\to \{g_i\}_{i\in |I|}$
be a morphism of
$\oplaxarrowsl{\scriptz}_a
\simeq \prod_{i\in |I|}
\oplaxarrowsl{\scriptz}_{a_i}$.
Then we have an oplax square $\sigma_i$:
\[ \xymatrix{
     x_i\ar[d]_{f_i}\ar[r]^{m_i^L}
     & y_i\ar[d]^{g_i}\\
     x_i'\ar@{=>}[ur]_{u_i}\ar[r]_{n_i} & y_i'\\
   }\]
for each $i\in |I|$,
where $m_i^L$ is a left adjoint $1$-morphism
in $\scriptz_{a_i}$.
By Lemma~\ref{lemma:oplax-arrow-L-cocart},
the morphism
$\{f_i\}_{i\in |I|}\to \{g_i\}_{i\in |I|}$
is $(s,t)_{sL,a}^{\otimes}$-coCartesian if and only if
the map 
$v_i: n_if_im_i^R\Longrightarrow g_i$
is invertible for each $i$,
where $v_i$ is a $2$-morphism
corresponding to $u_i$ under the equivalence
in Lemma~\ref{lemma:equivalence-mate-space}.
We can identify
the morphism
$\phi_*(\{f_i\}_{i\in |I|}\to \{g_i\}_{i\in |I|})$
with an oplax square $\sigma$:
\[ \xymatrix{
     \otimes_{\phi}x_i\ar[d]_{\otimes_{\phi}f_i}
     \ar[r]^{\otimes_{\phi}m_i^L}
     & \otimes_{\phi}y_i\ar[d]^{\otimes_{\phi}g_i}\\
     \otimes_{\phi}x_i'\ar[r]_{\otimes_{\phi}n_i}
     \ar@{=>}[ur]_{\otimes_{\phi}u_i}
     & \otimes_{\phi}y_i'.\\
   }\]
By Lemma~\ref{lemma:oplax-arrow-L-cocart},
$\sigma$ is an $(s,t)_{sL,b}^{\otimes}$-coCartesian morphism
if and only if
$v: (\otimes_{\phi}n_i)(\otimes_{\phi}f_i)
   (\otimes_{\phi}m_i)^R\Rightarrow
   \otimes_{\phi}g_i$
is invertible,
where $(\otimes_{\phi}m_i)^R$
is a right adjoint to
$\otimes_{\phi}m_i^L$.
Since there is an equivalence
$(\otimes_{\phi}m_i)^R\simeq
\otimes_{\phi}(m_i^R)$,
we see that $v\simeq \otimes_{\phi}v_i$
is invertible
if $v_i$ is invertible for each $i$.
\qed

\part{Convolution product
in monoidal $(\infty,2)$-categories}
\label{part:convolution-product}

In ordinary category theory
the endomorphism category of a map monoidale
has the structure of a duoidal category.
One of the two monoidal products
is given by composition of morphisms,
and the other is by convolution product. 

Let $\scriptz$ be an $\mathcal{O}$-monoidal
$(\infty,2)$-category
and let $A$ be a map $\mathcal{O}$-monoidale in $\scriptz$.
In \S\ref{section:map-duoidal}
we showed that
$\scriptz(A,A)^{\otimes}_{\bullet}$
corresponds to
a coCartesian $(\Delta^{\rm op},\mathcal{O})$-duoidal
$\infty$-category.
The $\Delta^{\rm op}$-monoidal
structure is given by composition
of morphisms in $\scriptz$.
In Part~\ref{part:convolution-product}
we show that the $\mathcal{O}$-monoidal
structure on
$\scriptz(A,A)^{\otimes}_{\bullet}$
is given by convolution product.

In \S\ref{section:convolution_product}
we introduce convolution product
on the mapping $\infty$-category $\scriptz(C,A)$
from an $\scripto$-comonoidale $C$
to an $\scripto$-monoidale $A$
in an $\scripto$-monoidal
$(\infty,2)$-category $\scriptz$.
We construct an $\mathcal{O}$-monoidal
$\infty$-category $\scriptz(C,A)^{\otimes}_{\rm conv}$
whose underlying $\infty$-category
is $\scriptz(C,A)$.
In \S\ref{section:duoidal-convolution}
we prove an equivalence
of $\scripto$-monoidal $\infty$-categories
between $\scriptz(A,A)^{\otimes}_1$
and $\scriptz(A^*,A)^{\otimes}_{\rm conv}$,
where
$A$ is a left adjoint $\scripto$-monoidale and
$A^*$ is the right adjoint $\scripto$-comonoidale
associated to $A$.

In order to prove this,
we need to establish
an equivalence
between some wide subcategories of
the twisted arrow $\infty$-category
and of the oplax arrow $\infty$-category.
For this purpose,
in \S\ref{section:twisted-squares}
we introduce a twisted square $\infty$-category,
which is a generalization of twisted arrow
$\infty$-category.
In \S\ref{section:twl-arrow-equivalence},
we show that there is such an equivalence 
by constructing a perfect pairing 
between them
by taking a suitable subcategory of
the twisted square $\infty$-category.

\section{Convolution product}
\label{section:convolution_product}

Let $\scriptz$ be an $\mathcal{O}$-monoidal
$(\infty,2)$-category
over a perfect operator category.
We take an $\mathcal{O}$-monoidale $A$
and an $\mathcal{O}$-comonoidale $C$ in $\scriptz$. 
In this section
we introduce a convolution product
on the mapping $\infty$-category $\scriptz(C,A)$,
which gives $\scriptz(C,A)$ the structure
of an $\mathcal{O}$-monoidal $(\infty,1)$-category.

\begin{definition}
[{\rm cf.~\cite[Definition~7.15]{HHLN2}}]\rm
For $n\ge 0$,
we denote by $[n]^{\rm op}\star_{\rm opl}[n]$
a strict $2$-category
informally depicted as \[\xymatrix{
     \overline{0}\ar[d]
     & \overline{1} \ar[l]\ar[d]\ar@{}[ld]|{\mbox{$\Longrightarrow$}}
     & \overline{2} \ar[l]\ar[d]\ar@{}[ld]|{\mbox{$\Longrightarrow$}}
     & \cdots \ar[l]\ar@{}[ld]|{\mbox{$\Longrightarrow$}}
     & \overline{n} \ar[l]\ar[d]
     \ar@{}[ld]|{\mbox{$\Longrightarrow$}} \\
     0 \ar[r]& 1 \ar[r]&
     2\ar[r]  & \cdots \ar[r]& n. }
\]
The functor $[n]\mapsto [n]^{\rm op}\star_{\rm opl}[n]$
forms a cosimplicial object in the category of
strict $2$-categories.
\end{definition}

\begin{definition}\rm
For an $(\infty,2)$-category $\scriptx$,
we define a simplicial $\infty$-groupoid
${\rm Tw}^l(\scriptx)_{\bullet}$ by
\[ {\rm Tw}^l(\scriptx)_n=
   {\rm Map}_{\bicat}([n]^{\rm op}\star_{\rm opl}[n],
   \scriptx). \]
As in the proof of \cite[Corollary~7.17]{HHLN2},
we see that
${\rm Tw}^l(\scriptx)_{\bullet}$ is
a complete Segal $\infty$-groupoid.
We denote by
${\rm Tw}^l(\scriptx)$
the $(\infty,1)$-category corresponding to
${\rm Tw}^l(\scriptx)_{\bullet}$.
\end{definition}

The inclusion functors
$[n]^{\rm op}\hookrightarrow [n]^{\rm op}\star_{\rm opl}[n]$
and
$[n]\hookrightarrow [n]^{\rm op}\star_{\rm opl}[n]$
induce a functor
\[ (s,t): {\rm Tw}^l(\scriptx)\longrightarrow
   u_1\scriptx^{\rm op}\times u_1\scriptx .\]
We understand that
an object of ${\rm Tw}^l(\scriptx)$
is a $1$-morphism $f: x\to y$ in $\scriptx$
and that a morphism of ${\rm Tw}^l(\scriptx)$
from $f: x\to y$ to $f':x'\to y'$ is a diagram $\sigma$:
\[ \begin{array}{ccc}
  x & \stackrel{g}{\longleftarrow} & x'\\
  \mbox{$\scriptstyle f$}\bigg\downarrow
  \phantom{\mbox{$\scriptstyle f$}}
  & \stackrel{\alpha}{\Longrightarrow} &
  \phantom{\mbox{$\scriptstyle f'$}}\bigg\downarrow
  \mbox{$\scriptstyle f'$}\\
  y & \stackrel{h}{\longrightarrow} & y'\\  
    \end{array}\] 
in $\scriptx$.
By \cite{Garcia-Stern}
(see also \cite[\S7]{HHLN2}),
the map $(s,t): {\rm Tw}^l(\scriptx)\longrightarrow
u_1\scriptx^{\rm op}\times u_1\scriptx$
is a coCartesian fibration
classified by the restricted mapping
$\infty$-category functor of $\scriptx$:
\[ {\scriptx}(-,-):
   u_1\scriptx^{\rm op}\times u_1\scriptx
   \longrightarrow \cat.\]
We note that a morphism $\sigma: f\to f'$
in ${\rm Tw}^l(\scriptx)$ is
$(s,t)$-coCartesian if and only if
the $2$-morphism $\alpha$ is invertible in $\scriptx$.

Now,
we let $\scriptz$
be an $\mathcal{O}$-monoidal $(\infty,2)$-category,
where $\mathcal{O}^{\otimes}$ is an $\infty$-operad
over a perfect operator category.
First,
we show that
${\rm Tw}^l(\scriptz)$ is an $\mathcal{O}$-monoidal
$(\infty,1)$-category.
The construction
$\scriptx\mapsto {\rm Tw}^l(\scriptx)$
determines a functor
${\rm Tw}^l: \bicat\to\cat$.
Since the functor ${\rm Tw}^l$
preserves finite products,
it induces a functor
${\rm Tw}^l: {\rm Mon}_{\mathcal{O}}(\bicat)\to
{\rm Mon}_{\mathcal{O}}(\cat)$.
Hence,
we obtain the following lemma.

\begin{lemma}
If $\scriptz$ is an $\mathcal{O}$-monoidal
$(\infty,2)$-category,
then ${\rm Tw}^l(\scriptz)$ has
an $\mathcal{O}$-monoidal structure.
\end{lemma}

We have a natural transformation
$(s,t): {\rm Tw}^l(-)\to u_1(-)^{\rm op}\times u_1(-)$.
Since the functors ${\rm Tw}^l(-)$
and $u_1(-)^{\rm op}\times u_1(-)$
preserve finite products,
we obtain the following lemma.

\begin{lemma}\label{lemma:Tw-strong-monoidal}
The map $(s,t): {\rm Tw}^l(\scriptz)\to
(u_1\scriptz)^{\rm op}\times u_1\scriptz$
is a strong $\mathcal{O}$-monoidal functor.
\end{lemma}

By Lemma~\ref{lemma:Tw-strong-monoidal},
we have a commutative diagram
in the $\infty$-category of $\infty$-operads
\begin{equation}\label{eq:o-monoida-triangle}
  \xymatrix{
  {\rm Tw}^l(\scriptz)^{\otimes}
  \ar[rr]^-{(s,t)^{\otimes}}\ar[dr]_p
  &&
  (u_1\scriptz^{\rm op})^{\otimes}\times
  u_1\scriptz^{\otimes}\ar[dl]^q\\
  & \mathcal{O}^{\otimes}, &}
  \end{equation}
where $p$ and $q$ are coCartesian fibrations,
and the top horizontal arrow $(s,t)^{\otimes}$
preserves coCartesian morphisms.
Next,
we show that
the map $(s,t)^{\otimes}:
{\rm Tw}^l(\scriptz)^{\otimes}\to
(u_1\scriptz^{\rm op})^{\otimes}\times u_1\scriptz^{\otimes}$
is a coCartesian fibration.

\begin{lemma}\label{lemma:Tw-u1-op-u1-cocart}
The map $(s,t)^{\otimes}:
{\rm Tw}^l(\scriptz)^{\otimes}\to
(u_1\scriptz^{\rm op})^{\otimes}\times u_1\scriptz^{\otimes}$
is a coCartesian fibration.
\end{lemma}

\proof
First,
we show that
$(s,t)^{\otimes}$
is a locally coCartesian fibration.
The maps $p$ and $q$ are coCartesian fibrations
and the map $(s,t)^{\otimes}$ carries
$p$-coCartesian morphisms to $q$-coCartesian morphisms.
Since the map
$(s,t): {\rm Tw}^l(\scriptz)\to
u_1\scriptz^{\rm op}\times u_1\scriptz$
is a coCartesian fibration,
we see that  
the map 
${\rm Tw}^l(\scriptz)^{\otimes}_{a}\to
(u_1\scriptz^{\rm op})^{\otimes}_a\times u_1\scriptz^{\otimes}_a$
induced on fibers 
is also a coCartesian fibration
for any $a\in\mathcal{O}^{\otimes}$.
Thus,
$(s,t)^{\otimes}$
is a locally coCartesian fibration
by the dual of \cite[Proposition~2.4.2.11]{Lurie1}.

Next,
we show that
$(s,t)^{\otimes}$
is in fact a coCartesian fibration.
By the dual of \cite[Proposition~2.4.2.8]{Lurie1},
it suffices to show that 
the induced functor
$\phi_*: {\rm Tw}^l(\scriptz)^{\otimes}_a\to
{\rm Tw}^l(\scriptz)^{\otimes}_b$
sends $(s,t)^{\otimes}_a$-coCartesian morphisms
to $(s,t)^{\otimes}_{b}$-coCartesian morphisms
for any morphism $\phi:a\to b$ in $\mathcal{O}^{\otimes}$.
We may assume that
$\phi: a\to b$ is active,
$p(a)=I$,
$a\simeq (a_i)_{i\in |I|}$,
and $a_i, b\in \mathcal{O}$.
We can identify a morphism $\sigma$
of ${\rm Tw}^l(\scriptz)^{\otimes}_a$
with a family $(\sigma_i)_{i\in |I|}$,
where $\sigma_i$ is a diagram
\[ \begin{array}{ccc}
  x_i & \stackrel{g_i}{\longleftarrow} & x'_i\\
  \mbox{$\scriptstyle f_i$}\bigg\downarrow
  \phantom{\mbox{$\scriptstyle f_i$}}
  & \stackrel{\alpha_i}{\Longrightarrow} &
  \phantom{\mbox{$\scriptstyle f'_i$}}\bigg\downarrow
  \mbox{$\scriptstyle f'_i$}\\
  y_i & \stackrel{h_i}{\longrightarrow} & y'_i\\  
    \end{array}\] 
in $\scriptz^{\otimes}_{a_i}$.
We see that the morphism 
$\phi_*\sigma\simeq
\otimes_{\phi}\sigma_i$ is a diagram
\[ \begin{array}{ccc}
  \otimes_{\phi}x_i & \stackrel{\otimes_{\phi}g_i}{\longleftarrow} &
  \otimes_{\phi}x'_i\\[1mm]  
  \mbox{$\scriptstyle \otimes_{\phi}f_i$}\bigg\downarrow
  \phantom{\mbox{$\scriptstyle \otimes_{\phi}f_i$}}
  & \stackrel{\otimes_{\phi}\alpha_i}{\Longrightarrow} &
  \phantom{\mbox{$\scriptstyle \otimes_{\phi}f'_i$}}\bigg\downarrow
  \mbox{$\scriptstyle \otimes_{\phi}f'_i$}\\
  y_i & \stackrel{\otimes_{\phi}h_i}{\longrightarrow} &
  \otimes_{\phi}y'_i\\  
    \end{array}\] 
in $\scriptz^{\otimes}_{b}$.
If $\sigma$ is an $(s,t)^{\otimes}_a$-coCartesian morphism,
then $\alpha_i$ is invertible for any $i\in |I|$.
This implies that $\otimes_{\phi}\alpha_i$ is also invertible.
By the characterization of coCartesian morphisms
of twisted arrow $\infty$-categories,
we see that $\phi_*\sigma$ is an $(s,t)^{\otimes}_{b}$-coCartesian
morphism. 
\qed

\bigskip

Now,
we introduce a convolution product
on the mapping $\infty$-category $\scriptz(C,A)$,
where $C\in {\rm coAlg}_{\mathcal{O}}(u_1\scriptz)$
is an $\mathcal{O}$-comonoidale
and $A\in {\rm Alg}_{\mathcal{O}}(u_1\scriptz)$
is an $\mathcal{O}$-monoidale
in $\scriptz$.
The objects $C$ and $A$ determine
morphisms
$C: \mathcal{O}^{\otimes} \to (u_1\scriptz^{\rm op})^{\otimes}$
and
$A: \mathcal{O}^{\otimes}\to u_1\scriptz^{\otimes}$
of $\infty$-operads over $\mathcal{O}^{\otimes}$,
respectively.
We consider the following pullback diagram
\begin{equation}\label{eq:defining-diagram-monoidal-Z-C-A}
  \begin{array}{ccc}
   \scriptz(C,A)^{\otimes}_{\rm conv}
  & \hbox to 10mm{\rightarrowfill} & {\rm Tw}^l(\scriptz)^{\otimes}\\
  \bigg\downarrow & &
  \phantom{\mbox{$\scriptstyle (s,t)^{\otimes}$}}
  \bigg\downarrow
  \mbox{$\scriptstyle (s,t)^{\otimes}$}\\
  \mathcal{O}^{\otimes}&
  \stackrel{(C,A)}{\hbox to 10mm{\rightarrowfill}} &
  (u_1\scriptz^{\rm op})^{\otimes}\times
  u_1\scriptz^{\otimes}.
  \end{array}
\end{equation} 
We can regard (\ref{eq:defining-diagram-monoidal-Z-C-A})
as a pullback diagram
in the $\infty$-category of $\infty$-operads over $\mathcal{O}^{\otimes}$.
Hence
the left vertical arrow
$\scriptz(C,A)^{\otimes}_{\rm conv}\to
\mathcal{O}^{\otimes}$ is a map of $\infty$-operads.

\begin{theorem}\label{thm:convolution-product}
Let $\scriptz$ be an $\mathcal{O}$-monoidal
$(\infty,2)$-category.
For an $\mathcal{O}$-monoidale
$A\in {\rm Alg}_{\mathcal{O}}(u_1\scriptz)$
and an $\mathcal{O}$-comonoidale
$C\in {\rm coAlg}_{\mathcal{O}}(u_1\scriptz)$
in $\scriptz$,
the mapping $(\infty,1)$-category $\scriptz(C,A)$
has an $\mathcal{O}$-monoidal structure.
\end{theorem}

\proof
It suffices to show that
the map
$\scriptz(C,A)^{\otimes}_{\rm conv}\to\mathcal{O}^{\otimes}$
is a coCartesian fibration.
Since 
the right vertical arrow in (\ref{eq:defining-diagram-monoidal-Z-C-A})
is a coCartesian fibration
by Lemma~\ref{lemma:Tw-u1-op-u1-cocart}, 
so is the left vertical arrow.
\qed

\begin{definition}\rm
We call the $\mathcal{O}$-monoidal
structure on $\scriptz(C,A)$
in Theorem~\ref{thm:convolution-product}
a convolution product.
\end{definition}

\begin{remark}\rm
We notice that
the construction of convolution product
is functorial.
We consider a functor
${\rm Mon}_{\mathcal{O}}(\bicat)\to \cat$
which assigns to an $\mathcal{O}$-monoidal
$(\infty,2)$-category $\scriptz$
the $(\infty,1)$-category 
${\rm coAlg}_{\scripto}(u_1\scriptz)^{\rm op}\times
{\rm Alg}_{\scripto}(u_1\scriptz)$.
We denote by
${\rm coAlg}_{\scripto}(u_1(-))^{\rm op}\times
{\rm Alg}_{\scripto}(u_1(-))\to
{\rm Mon}_{\scripto}(\bicat)$
a coCartesian fibration
obtained by unstraightening
the functor.
Informally speaking,
an object of ${\rm coAlg}_{\scripto}(u_1(-))^{\rm op}\times
{\rm Alg}_{\scripto}(u_1(-))$
is a triple $(\scriptz, C, A)$
in which $\scriptz$ is an $\scripto$-monoidal
$(\infty,2)$-category,
$C$ is an $\scripto$-comonoidale,
and $A$ is an $\scripto$-monoidale
in $\scriptz$.
We can construct a functor
\[ {\rm coAlg}_{\scripto}(u_1(-)^{\rm op})\times
   {\rm Alg}_{\scripto}(u_1(-))
   \longrightarrow
   {\rm Mon}_{\mathcal{O}}(\cat) \]
which assigns to $(\scriptz,C,A)$
the mapping $\infty$-category $\scriptz(C,A)$
equipped with the convolution $\scripto$-monoidal
structure.
\end{remark}

\section{Twisted square $\infty$-categories}
\label{section:twisted-squares}

In this section
we introduce a twisted square $\infty$-category
of an $(\infty,2)$-category,
which is a generalization of
twisted arrow $\infty$-category.
We construct a simplicial space 
$\TSop(X)_{\bullet}$ for an $\infty$-bicategory $X$,
which is equipped with a map
to $\twr{X}_{\bullet}\times\laxarrow{X}_{\bullet}$.
We defer a proof of $\TSop(X)_{\bullet}$
being in fact a complete Segal space
to \S\ref{section:proof-twisted-sq-complete-Segal}
below.

First,
in order to introduce a twisted square
$\infty$-category,
we define a cosimplicial object
$\ts^{\bullet}$
of scaled simplicial sets.


\begin{definition}\rm
We define a cosimplicial object
$\tsscplus$
of scaled simplicial sets.
For $n\ge 0$,
the underlying simplicial set
of $\tsscplus$
is $\Delta^2\times \Delta^n$.
We depict it by the following grid:
\[ \xymatrix{
  000\ar[r]\ar[d]
  & 001\ar[r]\ar[d]
  & \cdots\ar[r]
  & 00n\ar[d] \\
    010\ar[r]\ar[d]
  & 011\ar[r]\ar[d]
  & \cdots\ar[r]
  & 01n\ar[d] \\
    110\ar[r]
  & 111\ar[r]
  & \cdots\ar[r]
  & 11n. \\
}   \]
The set $T_{\tsscnplus}$
of thin $2$-simplices is given by
\[ \begin{array}{rccl}
  T_{\tsscnplus}
  &=& &
    \{\Delta^{\{ijk,ijk',ijk''\}} |\
    (i,j)=(0,0),(0,1),(1,1),\
    0\le k< k'< k'' \le n\}\\[2mm]
  & & \cup &
    \{\Delta^{\{00k,01k',01k''\}} |\
    0\le k\le  k'< k'' \le n\}\\[2mm]
  & & \cup &
    \{\Delta^{\{01k,01k',11k''\}} |\
      0\le k < k'\le  k'' \le n\}\\[2mm]
  & & \cup &
    \{\Delta^{\{00k,01k',11k''\}} |\
     0\le k\le k'\le  k'' \le n\}\\[2mm]
  & & \cup &
    \{\mbox{\rm degenerate}\}. \\     
\end{array}\]
We can verify that
$\tsscplus=\{\tsscnplus\}_{n\ge 0}$
forms a cosimplicial object.
\end{definition}

\begin{definition}\rm
We define a cosimplicial object
$\tsscminus$
of scaled simplicial sets.
For $n\ge 0$,
we consider the simplicial set
$\Delta^n\star\Delta^{n,{\rm op}}\star \Delta^n$.
We depict it by the following grid:
\[ \xymatrix{
  000\ar[r]\ar[d]
  & 001\ar[r]\ar[d]
  & \cdots\ar[r]
  & 00n\ar[d] \\
  100\ar[d]
  & 101\ar[l]\ar[d]
  & \cdots\ar[l]
  & 10n\ar[l]\ar[d] \\
   110\ar[r]
  & 111\ar[r]
  & \cdots\ar[r]
  & 11n. \\
}   \]
For a vertex $v=(i,k)$ of $\Delta^2\times \Delta^n$,
we set $\Omega(v)=00k$ if $i=0$,
$\Omega(v)=10k$ if $i=1$,
$\Omega(v)=11k$ if $i=2$.
For a simplex $\sigma$ of $\Delta^2\times\Delta^n$,
we denote by $\Omega(\sigma)$
the simplex of $\Delta^n\star\Delta^{n,{\rm op}}\star\Delta^n$
spanned by $\Omega(v)$
for vertices $v$ of $\sigma$.
For a subcomplex $K$ of $\Delta^2\times \Delta^n$,
we denote by $\Omega(K)$
the subcomplex of $\Delta^n\star\Delta^{n,{\rm op}}\star\Delta^n$
spanned by $\Omega(\sigma)$
for simplices $\sigma$ of $K$.

The underlying simplicial set
of $\tsscnminus$ 
is $\Omega(\Delta^2\times \Delta^n)$.
In other words,
it is spanned by
$(n+2)$-dimensional simplices
\[ \Delta^{\{000,\ldots,00k,10k,\ldots 10k',11k',\ldots,11n\}} \]
for $0\le k\le k'\le n$.
The set $T_{\tsscnminus}$
of thin $2$-simplices is given by
\[ \begin{array}{rccl}
  T_{\tsscnminus}
    &=& &
    \{\Delta^{\{ijk,ijk',ijk''\}}|\
    (i,j)=(0,0),(1,0),(1,1),\
    0\le k< k'< k'' \le n\}\\[2mm]
    & & \cup &
    \{\Delta^{\{00k,00k',10k''\}}|\
    0\le k<  k'\le k'' \le n\}\\[2mm]
  & & \cup &
    \{\Delta^{\{10k,11k',11k''\}}|\
      0\le k \le k'<  k'' \le n\}\\[2mm]
  & & \cup &
    \{\mbox{\rm degenerate}\}.     
\end{array}\]
We can verify that
$\tsscminus=\{\tsscnminus\}_{n \ge 0}$
forms a cosimplicial object.
\end{definition}

\begin{definition}\rm
We define a cosimplicial object
$\tsscbullet$
of scaled simplicial sets.
For $n\ge 0$,
we have monomorphisms
from $(\Delta^1\times\Delta^n)_{\flat}$
to $\tsscnplus$ and $\tsscnminus$, 
respectively,
given by
$(i,k)\mapsto iik$.
We define a scaled simplicial set
$\tsscn$
as a pushout
\[ \tsscn =
   \tsscnplus 
   \coprod_{(\Delta^1\times\Delta^n)_{\flat}}
   \tsscnminus.\]
We can verify that
$\tsscbullet =\{\tsscn \}_{n\ge 0}$
forms a cosimplicial object.
\end{definition}


\begin{notation}\rm
For a subcomplex $K$ of the underlying simplicial set
of $\tsscn$,
we denote by $K_{\dagger}$
the scaled simplicial set whose underlying simplicial set
is $K$ equipped with the induced scaling from
$\tsscn$.
\end{notation}

For scaled simplicial sets $A$ and $B$
with $B$ fibrant,
we recall that 
$\mapsc(A,B)$
is the underlying $\infty$-groupoid
of the $\infty$-bicategory ${\rm FUN}(A,B)$.

\begin{definition}\rm
For an $\infty$-bicategory $X$,
we define a simplicial space
$\TSop(X)_{\bullet}$
by 
\[ \TSop(X)_n =
   \mapsc(\tsscn, X).\]
\end{definition}

We show that
$\TSop(X)_{\bullet}$ is a complete Segal space
and hence it represents an $(\infty,1)$-category. 

\begin{proposition}\label{prop:TSZ-complete-Segal-space-orig}
For any $\infty$-bicategory $X$,
$\TSop(X)_{\bullet}$ is a complete Segal space.
\end{proposition}

We defer the proof of
Proposition~\ref{prop:TSZ-complete-Segal-space-orig}
to \S\ref{section:proof-twisted-sq-complete-Segal} below.

\begin{remark}\rm
By the assignment $X\mapsto
\TSop(X)_{\bullet}$,
we obtain a functor
$\TSop: (\setsc)^{\circ}\to
{\rm CSS}^{\circ}$,
where ${\rm CSS}$
is the category of simplicial spaces
equipped with the complete Segal space model structure.
We can verify that
$\TSop$ is an enriched functor between 
$(\setdelta)^{\circ}$-enriched categories.
By taking simplicial nerves,
we obtain a functor
${\rm TS}:\bicat\to\cat$
of $\infty$-categories.
\end{remark}


By restrictions to faces of $\tsscbullet$,
we construct maps out of $\TSop(X)_{\bullet}$.
For this purpose,
we set
\[ \begin{array}{rcl}
    \partial^{\rm T}\tsscn
    &=&(\Delta^{\{0,1\}}\times\Delta^n)_{\dagger},\\[2mm]
    \partial^{\rm F}\tsscn
    &=& (\Delta^{\{1,2\}}\times\Delta^n)_{\dagger},\\[2mm]
    \partial^{\rm R}\tsscn
    &=& \Omega(\Delta^{\{0,1\}}\times\Delta^n)_{\dagger},\\[2mm]
    \partial^{\rm B}\tsscn
    &=&\Omega(\Delta^{\{1,2\}}\times\Delta^n)_{\dagger}.\\[2mm]
\end{array}  \]
We can verify that
$\partial^{\rm F}\tsscbullet$,
$\partial^{\rm R}\tsscbullet$,
$\partial^{\rm T}\tsscbullet$,
$\partial^{\rm B}\tsscbullet$
form cosimplicial objects.

\begin{definition}\rm
For $f={\rm T}, {\rm F}, {\rm R}, {\rm B}$,
we define simplicial spaces
$\partial^f\TSop(X)_{\bullet}$ by
\[ \partial^f\TSop(X)_{\bullet}=
   \mapsc(\partial^f\tsscbullet,X). \]
\end{definition}

\begin{remark}\rm
For each $n\ge 0$,
there are natural isomorphisms
$\partial^{\rm T}\TSop(X)_n\cong
{\rm Fun}^{\rm opl}([n],{X})^{\simeq}$
and 
$\partial^{\rm F}\TSop(X)_n\cong
{\rm Fun}^{\rm lax}([n],X)^{\simeq}$
of simplicial sets.
Thus,
the simplicial spaces
$\partial^{\rm T}\TSop(X)_{\bullet}$
and 
$\partial^{\rm F}\TSop(X)_{\bullet}$
are complete Segal spaces.
\end{remark}

As for $\partial^{\rm R}\TSop(X)_{\bullet}$
and $\partial^{\rm B}\TSop(X)_{\bullet}$,
we show that they are also complete
Segal spaces by comparing with
twisted arrow $\infty$-categories.
We recall that
${\rm rev}: \Delta^{\rm op}\to\Delta^{\rm op}$
is a functor given by assigning to a nonempty finite ordered
set its reverse ordered set. 
For a simplicial object $S_{\bullet}$,
we denote by $S_{\bullet}^{\rm rev}$
a simplicial object given by
the composite functor $S_{\bullet}\circ {\rm rev}$.
We notice that
$\partial^{\rm B}\TSop(X)_{\bullet}$
is isomorphic to 
$\partial^{\rm R}\TSop(X)_{\bullet}^{\rm rev}$
as simplicial spaces.
Thus, it suffices to show that
$\partial^{\rm R}\TSop(X)_{\bullet}$
is a complete Segal space.

We recall the twisted arrow $\infty$-category
of an $(\infty,2)$-category.
Let $Q(\bullet)$ be a cosimplicial scaled simplicial set,
in which $Q(n)$ is given as in \cite[Definition~2.2]{Garcia-Stern}.
For a $\infty$-bicategory $X$,
we define a simplicial space
\[ \twr{X}_{\bullet}=\mapsc(Q(\bullet), X),\]
which is a complete Segal space
by \cite[Theorem~4.9]{Torii3}.

We have the following lemma.

\begin{lemma}
The simplicial spaces
$\partial^{\rm R}\TSop(X)_{\bullet}$
and
$\partial^{\rm B}\TSop(X)_{\bullet}$
are complete Segal spaces. 
There are equivalences
of complete Segal spaces
$\partial^{\rm R}\TSop(X)_{\bullet}
\simeq \twr{X}_{\bullet}$
and
$\partial^{\rm B}\TSop(X)_{\bullet}
\simeq \twr{X}_{\bullet}^{\rm rev}$.
\end{lemma}

\proof
Since we have an isomorphism
$\partial^{\rm B}\TSop(X)_{\bullet}\cong
\partial^{\rm R}\TSop(X)_{\bullet}^{\rm rev}$
of simplicial spaces,
it suffices to show that
$\partial^{\rm R}\TSop(X)_{\bullet}$
is a complete Segal space and
that there is an equivalence
of complete Segal spaces
between $\partial^{\rm R}\TSop(X)_{\bullet}$
and $\twr{X}_{\bullet}$.
We notice that
the cosimplicial scaled simplicial set
$\partial^{\rm R}\tsscbullet$ is isomorphic
to $T(\bullet)$ defined in \cite[Definition~5.1]{Torii3}.
The claim follows from
\cite[Proposition~5.8 and Theorem~5.9]{Torii3}.
\qed

\begin{definition}\label{def:partial-maps}
\rm
For $f={\rm T}, {\rm F}, {\rm R}, {\rm B}$,
we define maps of complete Segal spaces
by restrictions
\[   \partial^f: \TSop(X)_{\bullet}\longrightarrow
     \partial^f\TSop(X)_{\bullet}.\]
\end{definition}

\section{An equivalence between
$\twlsr{\scriptx}$ and $\oplaxarrowsl{\scriptx}$}
\label{section:twl-arrow-equivalence}

For an $(\infty,2)$-category $\scriptx$,
we set
\[ \twlsr{\scriptx}=
    \twl{\scriptx}\times_{(u_1\scriptx\times u_1\scriptx)}
    (u_1\scriptx^R\times u_1\scriptx).\]
We recall that
$\oplaxarrowsl{\scriptx}=
    \oplaxarrow{\scriptx}\times_{(u_1\scriptx\times u_1\scriptx)}
    (u_1\scriptx^L\times u_1\scriptx)$.
We would like to show that
there is an equivalence of $(\infty,1)$-categories
between $\twlsr{\scriptx}$
and $\oplaxarrowsl{\scriptx}$,
which is compatible with the dual equivalence
of $u_1\scriptx^R$ and $u_1\scriptx^L$
in Corollary~\ref{cor:dual-eq-u1-Z-L-R}.
We set 
\[ \twrsl{\scriptx}
   \simeq\twlsr{\scriptx^{\twoop}}^{\rm op},\qquad
   \laxarrowsr{\scriptx}
\simeq \oplaxarrowsl{\scriptx^{\twoop}}. \]
In this section
we show that 
there is a dual equivalence of $(\infty,1)$-categories
between $\twrsl{\scriptx}$ 
and $\laxarrowsr{\scriptx}$.
For this purpose,
we construct a perfect pairing between them
by taking a suitable subcategory of
the $(\infty,1)$-category of twisted squares
in $\scriptx$.


First, by using Proposition~\ref{prop:TSZ-complete-Segal-space-orig},
we define an $(\infty,1)$-category
of twisted squares in an $(\infty,2)$-category
$\scriptx$.

\begin{definition}\rm
Let $\scriptx$ be an $(\infty,2)$-category
represented by an $\infty$-bicategory $X$.
We denote by $\TSop(\scriptx)$
the $(\infty,1)$-category corresponding
to the complete Segal space $\TSop(X)_{\bullet}$.
We call $\TSop(\scriptx)$
the $(\infty,1)$-category of twisted squares in $\scriptx$.
\end{definition}

Informally speaking,
the objects of ${\rm TS}(\scriptx)$
are lax squares
\begin{equation}\label{eq:object-oplax-square-orig}
  \xymatrix{
      x(00) \ar[r]^{a_{00}}\ar[d]_{b_{00}} & 
      x(01)\ar[d]^{b_{01}}\ar@{=>}[dl]\\
      x(10) \ar[r]_{a_{10}} & x(11),\\
  }
\end{equation}
and morphisms are
diagrams 
\begin{equation}\label{eq:twisted-cube-in-Z-orig}
  \vcenter{
    \xymatrix{
     x(000) \ar[rr]^{m_{001}}\ar[dd]_{b_{000}}\ar[dr]^{a_{000}} & &
     x(001)\ar@{..>}[dd]^(.70){b_{001}}\ar[dr]^{a_{001}} & \\
     & x(010)\ar[dd]_(.30){b_{010}} \ar[rr]^(.30){m_{011}} &
     & x(011) \ar[dd]^{b_{011}}\\     
     x(100) \ar[dr]_{a_{100}} &
     & x(101)\ar@{..>}[ll]_(.30){n_{101}}\ar@{..>}[dr]^{a_{101}} &  \\
     & x(110)\ar[rr]_{m_{111}} &  & x(111)\\
   }}\end{equation}
with $2$-morphisms
\[ \xymatrix{
      x(000)\ar[r]^{a_{000}}\ar[d]_{b_{000}} &
      x(010)\ar[d]^{b_{010}}\ar@{=>}[dl]^{S_0}\\
      x(100)\ar[r]_{a_{100}}
      & x(110),\\}\qquad
    \xymatrix{
      x(001)\ar[r]^{a_{001}}\ar[d]_{b_{001}} &
      x(011)\ar[d]^{b_{011}}\ar@{=>}[dl]^{S_1}\\
      x(101)\ar[r]_{a_{101}}
      & x(111),\\}\qquad
    \xymatrix{
      x(000)\ar[r]^{m_{001}}\ar[d]_(0.4){a_{000}}
    & x(001)\ar[d]^(0.4){a_{001}}\\
    x(010)\ar[r]_{m_{011}}\ar@{=>}[ur]_{T_1} &
    x(011),} \]
\[  \xymatrix{
    x(100)\ar[d]_(0.4){a_{100}}
    \ar@{}[dr]|{\mbox{$\stackrel{B_1}{\Longleftarrow}$}}
    & x(101)\ar[d]^(0.4){a_{101}}\ar[l]_{n_{101}}\\
    x(110)\ar[r]_{m_{111}}  &
    x(111),}\qquad
    \xymatrix{
    x(010)\ar[r]^{m_{011}}\ar[d]_(0.4){b_{010}}
    & x(011)\ar[d]^(0.4){b_{011}}\ar@{=>}[dl]^{F_1}\\
    x(110)\ar[r]_{m_{111}} &
    x(111),}\qquad
   \xymatrix{
    x(000)\ar[r]^{m_{001}}\ar[d]_(0.4){b_{000}}
    \ar@{}[dr]|{\mbox{$\stackrel{R_1}{\Longrightarrow}$}}
    & x(001)\ar[d]^(0.4){b_{001}}\\
      x(100)  &
      x(101)\ar[l]^{n_{101}}}
   \]
in $\scriptx$
equipped with an equivalence
\[ \alpha_{01}:
   (B_1b_{001}m_{001}) \cdot (S_1m_{001})
   \cdot (b_{011}T_1)
   \stackrel{\simeq}{\longrightarrow}
   (m_{111}a_{100}R_1)\cdot (m_{111}S_0)
   \cdot (F_1a_{000}) \]
in ${\rm Map}_{\scriptx(x(000),x(111))}(f_{01},g_{01})$,
where $f_{01}=b_{011}m_{011}a_{000}$
and $g_{01}=m_{111}a_{100}n_{101}b_{001}m_{001}$.

The maps
$\partial^f: \TSop(X)_{\bullet}
\to\partial^f\TSop(X)_{\bullet}$
of complete Segal spaces for $f={\rm R}, {\rm F}$
induce a map of $(\infty,1)$-categories
\[ (\partial^{\rm R},\partial^{\rm F}):
   \TSop(\scriptx)\longrightarrow
   \twr{\scriptx}
   \times \laxarrow{\scriptx},\]
where $\laxarrow{\scriptx}$
is the lax arrow $\infty$-category
represented by
${\rm Fun}^{\rm lax}(\Delta^1_{\sharp},X)$
(cf.~\cite[Remark~7.7]{HHLN2}).

Now,
we study conditions
in which diagram~(\ref{eq:twisted-cube-in-Z-orig})
is an $(\partial^{\rm R},\partial^{\rm F})$-coCartesian morphism.
For this purpose,
we consider a diagram 
\begin{equation}\label{eq:twisted-cube-in-Z-2-orig}
  \vcenter{
    \xymatrix{
     x(000) \ar[rr]^{m_{001}}\ar[dd]_{b_{000}}\ar[dr]^{a_{000}} & &
     x(001) \ar[rr]^{m_{002}}\ar@{..>}[dd]^(.70){b_{001}}\ar[dr]^{a_{001}}
     & & x(002) \ar@{..>}[dd]^(.70){b_{002}}\ar[dr]^{a_{002}} \\
     & x(010)\ar[dd]_(.30){b_{010}} \ar[rr]^(.30){m_{011}} &
     & x(011)\ar[dd]^(.70){b_{011}} \ar[rr]^(.30){m_{012}}
     & & x(012)\ar[dd]^{b_{012}}   \\     
     x(100) \ar[dr]_{a_{100}} &
     & x(101)\ar@{..>}[ll]_(.30){n_{101}}\ar@{..>}[dr]^{a_{101}}
     & & x(102)\ar@{..>}[ll]_(.30){n_{102}}\ar@{..>}[dr]^{a_{102}}  \\
     & x(110)\ar[rr]_{m_{111}} &  & x(111)\ar[rr]_{m_{112}}
     & & x(112)\\
   }}\end{equation}
with $2$-morphisms
\[ \xymatrix{
      x(000)\ar[r]^{a_{000}}\ar[d]_{b_{000}} &
      x(010)\ar[d]^{b_{010}}\ar@{=>}[dl]^{S_0}\\
      x(100)\ar[r]_{a_{100}}
      & x(110),\\}\qquad
    \xymatrix{
      x(001)\ar[r]^{a_{001}}\ar[d]_{b_{001}} &
      x(011)\ar[d]^{b_{011}}\ar@{=>}[dl]^{S_1}\\
      x(101)\ar[r]_{a_{101}}
      & x(111),\\}\qquad
    \xymatrix{
      x(002)\ar[r]^{a_{002}}\ar[d]_{b_{002}} &
      x(012)\ar[d]^{b_{012}}\ar@{=>}[dl]^{S_2}\\
      x(102)\ar[r]_{a_{102}}
      & x(112),\\}\]
\[ \xymatrix{
    x(000)\ar[r]^{m_{001}}\ar[d]_(0.4){a_{000}}
    & x(001)\ar[d]^(0.4){a_{001}}\ar[r]^{m_{002}}
    & x(002)\ar[d]^(0.4){a_{002}}\\
    x(010)\ar[r]_{m_{011}}\ar@{=>}[ur]_{T_1} &
    x(011)\ar[r]_{m_{012}}\ar@{=>}[ur]_{T_2} & x(012),}\qquad
   \xymatrix{
    x(100)\ar[d]_(0.4){a_{100}}
    \ar@{}[dr]|{\mbox{$\stackrel{B_1}{\Longleftarrow}$}}
    & x(101)\ar[d]^(0.4){a_{101}}\ar[l]_{n_{101}}
    \ar@{}[dr]|{\mbox{$\stackrel{B_2}{\Longleftarrow}$}}
    & x(102)\ar[d]^(0.4){a_{102}}\ar[l]_{n_{102}}\\
    x(110)\ar[r]_{m_{111}}  &
    x(111)\ar[r]_{m_{112}} & x(112),}
   \]
\[ \xymatrix{
    x(010)\ar[r]^{m_{011}}\ar[d]_(0.4){b_{010}}
    & x(011)\ar[d]^(0.4){b_{011}}\ar[r]^{m_{012}}\ar@{=>}[dl]^{F_1}
    & x(012)\ar[d]^(0.4){b_{012}}\ar@{=>}[dl]^{F_2}\\
    x(110)\ar[r]_{m_{111}} &
    x(111)\ar[r]_{m_{112}} & x(112),}\qquad
   \xymatrix{
    x(000)\ar[r]^{m_{001}}\ar[d]_(0.4){b_{000}}
    \ar@{}[dr]|{\mbox{$\stackrel{R_1}{\Longrightarrow}$}}
    & x(001)\ar[d]^(0.4){b_{001}}\ar[r]^{m_{002}}
    \ar@{}[dr]|{\mbox{$\stackrel{R_2}{\Longrightarrow}$}}
    & x(002)\ar[d]^(0.4){b_{002}}\\
    x(100)  &
    x(101)\ar[l]^{n_{101}} & x(102)\ar[l]^{n_{102}}.}
   \]
in $\scriptx$.
We set 
\[ \begin{array}{rcl}
  f_{01}&=&b_{011}m_{011}a_{000},\\
  f_{12}&=&b_{012}m_{012}a_{001},\\
  f_{02}&=&b_{012}m_{012}m_{011}a_{000}.\\
\end{array}\]
and
\[ \begin{array}{rcl}
  g_{01}&=&m_{111}a_{100}n_{101}b_{001}m_{001},\\
  g_{12}&=&m_{112}a_{101}n_{102}b_{002}m_{002}, \\
  g_{02}&=&m_{112}m_{111}a_{100}n_{101}n_{102}b_{002}m_{002}m_{001}.\\
  \end{array}\]
By composition with $T_1$ and $B_1$ ,
we obtain a map of spaces
\[ B_1(-)T_1: {\rm Map}_{\scriptx(x(001),x(112))}(f_{12},g_{12})
   \longrightarrow
       {\rm Map}_{\scriptx(x(000),x(112))}(f_{02},g_{02}).\]

\begin{lemma}\label{lemma:T-B-equivalence-orig}
We assume that $m_{001}$ admits a right adjoint
$m_{001}^R: x(001)\to x(000)$,
that a mate $m_{011}a_{000}m_{001}^R\Rightarrow a_{001}$
of $T_1$ is invertible
and that the $2$-morphism $B_1$ is invertible.
Then the map $B_1(-)T_1$ is an equivalence.
\end{lemma}

\proof
By the assumption that
$B_1$ is invertible,
the induced map
\[ (B_1)_*: \map{\scriptx(x(000),x(112))}(f_{02},g_{12}m_{001})
        \stackrel{\simeq}{\longrightarrow}
        \map{\scriptx(x(000),x(112))}(f_{02},g_{02})\]
is an equivalence.
Since $m_{001}$ admits a right adjoint $m_{001}^R$,
we have an equivalence
\[ \map{\scriptx(x(001),x(112))}(f_{02}m_{001}^R,g_{12})
   \simeq
   \map{\scriptx(x(000),x(112))}(f_{02},g_{12}m_{001})\]
by Lemma~\ref{lemma:equivalence-mate-space}.
The lemma follows from
the fact that
there is an equivalence $f_{12}\simeq f_{02}m_{001}^R$
since a mate of $T_1$ is invertible.
\qed

\bigskip

For simplicity,
we write
\[ \scriptx(f_{01},g_{01}), \quad
   \scriptx(f_{12},g_{12}), \quad
   \scriptx(f_{02},g_{02}) \]
for
${\rm Map}_{\scriptx(x(000),x(111))}(f_{01},g_{01})$,
${\rm Map}_{\scriptx(x(001),x(112))}(f_{12},g_{12})$,
${\rm Map}_{\scriptx(x(000),x(112))}(f_{02},g_{02})$,
respectively.
We set
\[ \begin{array}{rcl}
    \Phi_{01}&=&
    B_1b_{001}m_{001}\cdot S_1m_{001}\cdot b_{011}T_1\\
    \Phi_{12}&=&
    B_2b_{002}m_{002}\cdot S_2m_{002}\cdot b_{012}T_2\\
    \Phi_{02}&=&
    B_{12}b_{002}m_{002}m_{001}\cdot S_2m_{002}m_{001}\cdot
    b_{012}T_{12}\\
    \end{array}\]
and
\[ \begin{array}{rcl}
    \Psi_{01}&=& m_{111}a_{100}R_1\cdot m_{111}S_0\cdot F_1a_{000} \\
    \Psi_{12}&=& m_{112}a_{101}R_2\cdot m_{112}S_1\cdot F_2a_{001} \\
    \Psi_{02}&=& m_{112}m_{111}a_{100}R_{12}\cdot m_{112}m_{111}S_0
    \cdot F_{12}a_{000} \\
   \end{array}\]
where
\[ \begin{array}{rcl}
     B_{12}&=& m_{112}B_1n_{102}\cdot B_2,\\
     T_{12}&=& T_2m_{001}\cdot m_{012}T_1,\\
     R_{12}&=& n_{101}R_2m_{001} \cdot R_1,\\
     F_{12}&=& m_{112}F_1\cdot F_2 m_{011}.\\
   \end{array}\]  
Note that $\Phi_{01},\Psi_{01} \in \scriptx(f_{01},g_{01})$,
$\Phi_{12},\Psi_{12} \in \scriptx(f_{12},g_{12})$, and
$\Phi_{02},\Psi_{02} \in \scriptx(f_{02},g_{02})$.

Suppose that we have morphism~(\ref{eq:twisted-cube-in-Z-orig}).
This determines an equivalence
$\alpha_{01}: \Phi_{01}\stackrel{\simeq}{\to}
           \Psi_{01}$
in $\scriptx(f_{01},g_{01})$.
Using this,
we construct a map of spaces
\[ {\rm Map}_{\scriptx(f_{12},g_{12})}(\Phi_{12},\Psi_{12})
   \longrightarrow
   {\rm Map}_{\scriptx(f_{02},g_{02})}(\Phi_{02},\Psi_{02})\]
as follows:
An equivalence $\alpha_{12}:
\Phi_{12}\stackrel{\simeq}{\to} \Psi_{12}$
induces an equivalence
$\Phi_{02}= B_1 \cdot \Phi_{12}\cdot T_1
           \stackrel{\simeq}{\to}
           B_1\cdot \Psi_{12}\cdot T_1$.
We obtain an equivalence
$B_1\cdot \Psi_{12}\cdot T_1\stackrel{\simeq}{\to}
   R_2\cdot \Phi_{01}\cdot F_2$
by using the canonical equivalences
$B_1\cdot R_2\simeq R_2\cdot B_1$
and
$F_2\cdot T_1\simeq T_1\cdot F_2$.
Composing with the equivalence
$\alpha_{01}: \Phi_{01}\stackrel{\simeq}{\to}
\Psi_{01}$,
we obtain an equivalence
\[ \Phi_{02}= B_1 \cdot \Phi_{12}\cdot T_1
   \stackrel{\alpha_{12}}{\longrightarrow}
   B_1 \cdot \Psi_{12}\cdot T_1
   \stackrel{\simeq}{\longrightarrow}
   R_2 \cdot \Phi_{01}\cdot F_2
   \stackrel{\alpha_{01}}{\longrightarrow}
   R_2 \cdot \Psi_{01}\cdot F_2
   =\Psi_{02}.\]

\begin{lemma}\label{lemma:Ph-Psi-equivalence-orig}
The map
$\map{\scriptx(f_{12,}g_{12})}(\Phi_{12},\Psi_{12})
   \to
   \map{\scriptx(f_{02},g_{02})}(\Phi_{02},\Psi_{02})$
is an equivalence.
\end{lemma}

\proof
By Lemma~\ref{lemma:T-B-equivalence-orig},
there is an equivalence
\[ \map{\scriptx(f_{12},g_{12})}(\Phi_{12},\Psi_{12})
   \stackrel{\simeq}{\longrightarrow}
   \map{\scriptx(f_{02},g_{02})}
   (B_1 \cdot \Phi_{12}\cdot T_1,
    B_1 \cdot \Psi_{12}\cdot T_1).\]
By composing with the equivalence
$B_1 \cdot \Psi_{12}\cdot T_1
 \stackrel{\simeq}{\to}
 R_2 \cdot \Phi_{01}\cdot F_2
 \stackrel{\alpha_{01}}{\to}
 R_2 \cdot \Psi_{01}\cdot F_2
 =\Psi_{02}$,
we obtain the desired equivalence. 
\qed

\begin{proposition}
If the $1$-morphism $m_{001}$ admits a right adjoint
$m_{001}^R: x(001)\to x(000)$,
a mate $m_{011}a_{000}m_{001}^R \Rightarrow a_{001}$
of $T_1$ is invertible,
and the $2$-morphism $B_1$ is invertible,
then diagram~{\rm (\ref{eq:twisted-cube-in-Z-orig})}
is an $(\partial^{\rm R},\partial^{\rm F})$-coCartesian morphism.
\end{proposition}

\proof
For simplicity,
we set
$\TSop=\TSop(\scriptx)$,
$\tw=\twr{\scriptx}$
and $\arrow=\laxarrow{\scriptx}$.
We have to show that
the following commutative diagram
\[ \begin{array}{ccc}
  \map{\TSop}(S_1,S_2)
   &\longrightarrow&
  \map{\TSop}(S_0,S_2)\\
  \bigg\downarrow & & \bigg\downarrow \\
  \map{\tw\times\arrow}
  ((b_{001},b_{011}),(b_{002},b_{012}))
  &\longrightarrow&
  \map{\tw\times\arrow}
  ((b_{000},b_{010}),(b_{002},b_{012}))
  \\
    \end{array}   \]
is pullback.       
We take fibers of the vertical arrows
at $(R_2,F_2)$ and $(R_{12}, F_{12})$,
respectively,
where $R_{12}= n_{101}R_2m_{001} \cdot R_1$
and
$F_{12}= m_{112}F_1 \cdot F_2m_{011}$.
It suffices to show that
the induced map
\[ {\rm Map}_{{\rm TS}}(S_1,S_2)_{(R_2,F_2)}
   \longrightarrow
   {\rm Map}_{{\rm TS}}(S_0,S_2)_{(R_{12},F_{12})}\]
is an equivalence.

Taking top and bottom faces of
diagram~(\ref{eq:twisted-cube-in-Z-2-orig}),
we obtain a commutative diagram
\begin{equation}\label{eq:diagram-pullback-boundary-faces-orig}
  \begin{array}{ccc}
    \map{{\rm TS}}(S_1,S_2)_{(R_2,F_2)}
    &\longrightarrow&
    \map{{\rm TS}}(S_0,S_2)_{(R_{12},F_{12})}\\
    \bigg\downarrow & & \bigg\downarrow \\
    {\rm top}_{12}\times {\rm btm}_{12}
    & \longrightarrow &
    {\rm top}_{02}\times {\rm btm}_{02},\\
  \end{array}
\end{equation}
where
\[ \begin{array}{rcl}
  {\rm top}_{12}&=&
  \map{\scriptx(x(001),x(012))}(a_{002}m_{002},m_{012}a_{001}),\\
  {\rm top}_{02}&=&
  \map{\scriptx(x(000),x(012))}
  (a_{002}m_{002}m_{001},m_{012}m_{011}a_{000}),\\  
\end{array}\]
and
\[ \begin{array}{rcl} 
  {\rm btm}_{12}&=&
  \map{\scriptx(x(102),x(112))} 
    (m_{112}a_{101}n_{102},a_{102}),\\
  {\rm btm}_{02}&=&
  \map{\scriptx(x(102),x(112))}
  (m_{112}m_{111}a_{100}n_{101}n_{102},a_{102}).\\
   \end{array}\]
The map
${\rm top}_{12}\to {\rm top}_{02}$ is an equivalence
by the assumptions that
$m_{001}$ admits a right adjoint
and that a mate of $T_1$ is invertible.
The map
${\rm btm}_{12}\to {\rm btm}_{02}$ is also an equivalence
by the assumption that $B_1$ is invertible.
Thus,
it suffices to show that
diagram~(\ref{eq:diagram-pullback-boundary-faces-orig})
is pullback.

We take $(T_2,B_2)\in {\rm top}_{12}\times {\rm btm}_{12}$
and
$(T_{12},B_{12})\in {\rm top}_{02}\times {\rm btm}_{02}$,
where $T_{12}= T_2m_{001} \cdot m_{012}T_1$
and
$B_{12}= m_{112}B_1n_{102} \cdot B_2$.
It suffices to show that
the induced map on fibers is
an equivalence.
We can identify it
with the map
${\rm Map}_{\scriptx(f_{12,}g_{12})}(\Phi_{12},\Psi_{12})
   \to
   {\rm Map}_{\scriptx(f_{02},g_{02})}(\Phi_{02},\Psi_{02})$,
which is an equivalence by
Lemma~\ref{lemma:Ph-Psi-equivalence-orig}.
\qed

\bigskip

We define
$\TSopl(\scriptx)$ to be 
the wide subcategory of ${\rm TS}(\scriptx)$
spanned by
those morphisms~(\ref{eq:twisted-cube-in-Z-orig})
in which $m_{001}$ is a left adjoint $1$-morphism
in $\scriptx$.
Then,
we obtain a map
\[ (\partial^{\rm R},\partial^{\rm F})_{L}:
   \TSopl(\scriptx)\longrightarrow
   \twrsl{\scriptx}
   \times \laxarrow{\scriptx}\]
by restricting $(\partial^{\rm R},\partial^{\rm F})$
to $\TSopl(\scriptx)$.
We notice that there is a pullback diagram
\[ \begin{array}{ccc}
    \TSopl(\scriptx)
    & \longrightarrow &
    \TSop(\scriptx) \\
    \mbox{$\scriptstyle
      (\partial^{\rm R},\partial^{\rm F})_L$}
    \bigg\downarrow
    \phantom{\mbox{$\scriptstyle
        (\partial^{\rm R},\partial^{\rm F})_L$}}
    & &
    \phantom{\mbox{$\scriptstyle
        (\partial^{\rm R},\partial^{\rm F})$}}
    \bigg\downarrow
    \mbox{$\scriptstyle
      (\partial^{\rm R},\partial^{\rm F})$}\\
    \twrsl{\scriptx}\times
    \laxarrow{\scriptx}
    & \longrightarrow &
   \twr{\scriptx}\times
    \laxarrow{\scriptx}.\\    
   \end{array}\] 

\begin{corollary}\label{cor:twisted-st-R-coCart-orig}
The map $(\partial^{\rm R},\partial^{\rm F})_{L}:
    \TSopl(\scriptx)\to
   \twrsl{\scriptx}
   \times \laxarrow{\scriptx}$
is a coCartesian fibration.
Diagram~\mbox{\rm (\ref{eq:twisted-cube-in-Z-orig})}
is a $(\partial^{\rm R},\partial^{\rm F})_{L}$-coCartesian
morphism if and only if
a mate $m_{011}a_{000}m_{001}^R\Rightarrow a_{001}$
of $T_1$ and the $2$-morphism
$B_1: a_{101}\Rightarrow m_{111}a_{100}n_{101}$
are invertible.
\end{corollary}

We define 
$\TSoprl(\scriptx)$ by
the following pullback diagram
\[ \begin{array}{ccc}
    \TSoprl(\scriptx)
    & \longrightarrow &
    \TSopl(\scriptx) \\
    \mbox{$\scriptstyle
      (\partial^{\rm R},\partial^{\rm F})_{L,R}$}
    \bigg\downarrow
    \phantom{\mbox{$\scriptstyle
        (\partial^{\rm R},\partial^{\rm F})_{L,R}$}}
    & &
    \phantom{\mbox{$\scriptstyle
        (\partial^{\rm R},\partial^{\rm F})_L$}}
    \bigg\downarrow
    \mbox{$\scriptstyle
      (\partial^{\rm R},\partial^{\rm F})_L$}\\
    \twrsl{\scriptx}\times
    \laxarrowsr{\scriptx}
    & \longrightarrow &
   \twrsl{\scriptx}\times
    \laxarrow{\scriptx}.\\    
   \end{array}\] 
By Corollary~\ref{cor:twisted-st-R-coCart-orig},
the left vertical arrow
$(\partial^{\rm R},\partial^{\rm F})_{L,R}$
is a coCartesian fibration.
We denote by $\TSoprlcocart(\scriptx)$
the wide subcategory of $\TSoprl(\scriptx)$
spanned by
$(\partial^{\rm R},\partial^{\rm F})_{L,R}$-coCartesian
morphisms.
By the restriction of
$(\partial^{\rm R},\partial^{\rm F})_{L,R}$
to $\TSoprlcocart(\scriptx)$,
we obtain a left fibration
\begin{equation}\label{eq:pairing-coCart-orig}
  (\partial^{\rm R},\partial^{\rm F})_{L,R}^{\rm cocart}:
   \TSoprlcocart(\scriptx)\to
   \twrsl{\scriptx}
   \times \laxarrowsr{\scriptx}.
   \end{equation}
Thus,
we obtain a pairing between
$\twrsl{\scriptx}$
and $\laxarrowsr{\scriptx}$.

Pairing~(\ref{eq:pairing-coCart-orig})
is not a perfect pairing in general.
By restricting
$\TSoprlcocart(\scriptx)$
to a full subcategory,
we construct a perfect pairing
between 
$\twrsl{\scriptx}$
and $\laxarrowsr{\scriptx}$.
We denote by
$\TSoprlpair(\scriptx)$
the full subcategory of 
$\TSoprlcocart(\scriptx)$
spanned by those objects~(\ref{eq:object-oplax-square-orig})
in which the $1$-morphism $a_{00}: x(00)\to x(01)$
is right adjoint.
By restricting
$(\partial^{\rm R},\partial^{\rm F})_{R,L}^{\rm cocart}$ to
$\TSoprlpair(\scriptx)$,
we obtain a map
\[ (\partial^{\rm R},\partial^{\rm F})_{L,R}^{\rm pair}:
   \TSoprlpair(\scriptx)\to
   \twrsl{\scriptx}
   \times \laxarrowsr{\scriptx}.
\]

\begin{lemma}\label{lemma-s-t-pair-Ts-pairing-orig}
The map $(\partial^{\rm R},\partial^{\rm F})_{L,R}^{\rm pair}$
is a left fibration.
\end{lemma}

\proof
Since $(\partial^{\rm R},\partial^{\rm F})_{L,R}^{\rm coCart}$
is a left fibration,
it suffices to show that if
the source of a morphism
of $\TSoprlcocart(\scriptx)$
is in $\TSoprlpair(\scriptx)$,
then the target is also
in $\TSoprlpair(\scriptx)$.
We suppose that 
diagram ~(\ref{eq:twisted-cube-in-Z-orig})
is an
$(\partial^{\rm R},\partial^{\rm F})_{L,R}$-coCartesian
morphism
and
that $a_{000}$ is a right adjoint $1$-morphism.
Then
$a_{001}$ is also a right adjoint $1$-morphism
since $m_{001}$ is left adjoint
and
a mate $m_{001}a_{000}m_{001}^R \Rightarrow a_{001}$
of $T_1$ is invertible.
\qed

\begin{theorem}\label{thm:perfect-pairing-TS}
The map
$(\partial^{\rm R},\partial^{\rm F})_{L,R}^{\rm pair}:
\TSoprlpair(\scriptx)\to
\twrsl{\scriptx}\times
\laxarrowsr{\scriptx}$ is a perfect
pairing.
\end{theorem}

\proof
By Lemma~\ref{lemma-s-t-pair-Ts-pairing-orig},
the map $(\partial^{\rm R},\partial^{\rm F})_{L,R}^{\rm pair}$ is a pairing
of $\infty$-categories between
$\twrsl{\scriptx}$ and 
$\laxarrowsr{\scriptx}$.
In order to prove that it is perfect,
by \cite[Corollary~5.2.1.22]{Lurie2},
we need to show that
it is both left and right representable,
and that an object of
$\TSoprlpair(\scriptx)$
is left universal if and only if it is
right universal.

First,
we shall show that
$(\partial^{\rm R},\partial^{\rm F})_{L,R}^{\rm pair}$ is left representable.
We regard a $1$-morphism $b: x\to y$
in $\scriptx$ as an object of
$\twrsl{\scriptx}$.
We set
$\mathcal{C}(b)=\{b\}
\times_{\twrsl{\scriptx}}
\TSoprlpair(\scriptx)$.
We regard
${\rm id}_b$ as 
an object of $\mathcal{C}(b)$
and show that it is initial.
We take an object $\sigma$
of $\mathcal{C}(b)$ that is represented by
a lax square
\[  \xymatrix{
     x \ar[r]^{a}\ar[d]_{b}
     & x' \ar[d]^{b'}\ar@{=>}[dl]^{\sigma}\\
     y \ar[r]_{a'}
     & y'.\\
  }\]
A morphism from ${\rm id}_b$ to $\sigma$
in $\mathcal{C}(b)$ is represented by a cube
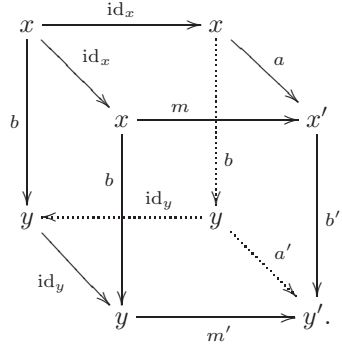
\[  \vcenter{
    \xymatrix{
     x \ar[rr]^{{\rm id}_x}\ar[dd]_{b}\ar[dr]^{{\rm id}_x} & &
     x\ar@{..>}[dd]^(.70){b}\ar[dr]^{a} & \\
     & x\ar[dd]_(.30){b} \ar[rr]^(.30){m} &
     & x' \ar[dd]^{b'}\\     
     y \ar[dr]_{{\rm id}_y} &
     & y\ar@{..>}[ll]_(.30){{\rm id}_y}\ar@{..>}[dr]^{a'} &  \\
     & y\ar[rr]_{m'} &  & y'.\\
}} \]
We note that 
the $2$-morphisms
$T_1: m\Rightarrow a$
and
$B_1: a'\Rightarrow m'$
are invertible.
By observing
that there is an equivalence
between
${\rm Map}_{\mathcal{C}(b)}({\rm id}_b,\sigma)$
and
${\rm Map}_{\scriptx(x,y')}(b'a,a'b)_{\sigma/}$,
we see that 
${\rm Map}_{\mathcal{C}(b)}({\rm id}_b,\sigma)$
is contractible.

Next,
we shall show that
$(\partial^{\rm R},\partial^{\rm F})_{L,R}^{\rm pair}$ is right representable.
We regard a $1$-morphism $b: x\to y$
in $\scriptx$ as an object of
$\laxarrowsr{\scriptx}$.
We set
$\mathcal{D}(b)=\{b\}
\times_{\laxarrowsr{\scriptx}}
\TSoprlpair(\scriptx)$.
We regard 
${\rm id}_b$ as an object of
$\mathcal{D}(b)$ and 
show that it is initial.
We take an object $\tau$
of $\mathcal{D}(b)$ that is represented by
a lax square
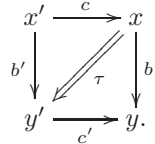
\[  \xymatrix{
     x' \ar[r]^{c}\ar[d]_{b'}
     & x \ar[d]^{b}\ar@{=>}[dl]^{\tau}\\
     y' \ar[r]_{c'}
     & y.\\
  }\]
A morphism from ${\rm id}_b$ to $\tau$
in $\mathcal{D}(b)$ is represented by a cube
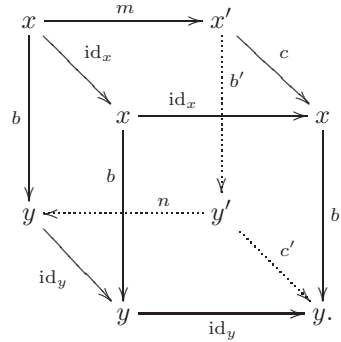
\[  \vcenter{
    \xymatrix{
     x \ar[rr]^{m}\ar[dd]_{b}\ar[dr]^{{\rm id}_x} & &
     x'\ar@{..>}[dd]^(.30){b'}\ar[dr]^{c_{}} & \\
     & x\ar[dd]_(.30){b} \ar[rr]^(.30){{\rm id}_x} &
     & x \ar[dd]^{b}\\     
     y \ar[dr]_{{\rm id}_y} &
     & y'\ar@{..>}[ll]_(.30){n}\ar@{..>}[dr]^{c'} &  \\
     & y\ar[rr]_{{\rm id}_y} &  & y.\\
}} \]
We note that 
a mate $m^R \Rightarrow c$ of $T_1$
and the $2$-morphism
$B_1: c'\Rightarrow n$
are invertible.
By observing
that there is an equivalence
between
${\rm Map}_{\mathcal{D}(b)}({\rm id}_b,\tau)$
and
${\rm Map}_{\scriptx(x',y)}(bc,c'b')_{/\tau}$,
we see that 
${\rm Map}_{\mathcal{D}(b)}({\rm id}_b,\tau)$
is contractible.

By the above argument,
we see that an object
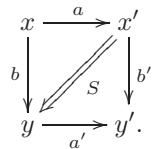
\[  \xymatrix{
     x \ar[r]^{a}\ar[d]_{b}
     & x'\ar[d]^{b'}\ar@{=>}[dl]^{S}\\
     y \ar[r]_{a'}
     & y'.\\
  }\]
of 
$\TSoprlpair(\scriptx)$
is left universal if and only
if $a$ and $a'$ are equivalence $1$-morphisms
and $S$ is an invertible $2$-morphism
in $\scriptx$
if and only if
it is right universal.
\qed

\bigskip

Theorem~\ref{thm:perfect-pairing-TS}
implies that there is a dual equivalence
of $(\infty,1)$-categories
between $\twrsl{\scriptx}$
and $\laxarrowsr{\scriptx}$.
Next,
we would like to show that
this equivalence is compatible with
the dual equivalence in Corollary~\ref{cor:dual-eq-u1-Z-L-R}.

Taking the top and bottom faces
in (\ref{eq:twisted-cube-in-Z-orig}),
we obtain a map
\[ (\partial^{\rm T},\partial^{\rm B}):
   \TSop(\scriptx)\longrightarrow
          \oplaxarrow{\scriptx}
          \times \twl{\scriptx^{\twoop}}.\]
By restriction,
we obtain a map of pairings
\begin{equation}\label{eq:map-pairing-orig}
  \begin{array}{ccc}
    \TSoprlpair(\scriptx)
    &\stackrel{(\partial^{\rm T},\partial^{\rm B})}{\longrightarrow}&
    {\rm Pair}(\scriptx)
    \times \twl{\scriptx^{\twoop}}\\[2mm]
    \mbox{$\scriptstyle (\partial^{\rm R},\partial^{\rm F})_{L,R}^{\rm pair}$}
    \bigg\downarrow
    \phantom{\mbox{$\scriptstyle (\partial^{\rm R},\partial^{\rm F})_{L,R}^{\rm pair}$}}
    & &
    \phantom{\mbox{$\scriptstyle (s,t)^{\rm pair}\times (s,t)$}}
    \bigg\downarrow
    \mbox{$\scriptstyle (s,t)^{\rm pair}\times (s,t)$}\\
    \twrsl{\scriptx}
    \times \laxarrowsr{\scriptx}
    &\longrightarrow &
    (u_1\scriptx^{L}\times u_1\scriptx^R)
    \times
    (u_1\scriptx^{\rm op}\times u_1\scriptx).\\
  \end{array}
\end{equation}

\begin{lemma}\label{lemma:map-pairing-representable}
The map of pairings in {\rm (\ref{eq:map-pairing-orig})}
is left and right representable
in the sense of {\rm \cite[Variant~5.2.1.16]{Lurie2}}.
\end{lemma}

\proof
Since the both pairings are perfect,
it suffices to show that
it is left representable.
From the description of left
representable objects in the pairings
in the proofs of 
Theorems~\ref{theorem:L-R-perfect-pairing}
and \ref{thm:perfect-pairing-TS},
we see that
the map $(\partial^{\rm T},\partial^{\rm B}): \TSoprlpair(\scriptx)
\to {\rm Pair}(\scriptx)
\times \twl{\scriptx^{\twoop}}$
carries left representable objects
to left representable objects.
\qed

\bigskip

Making use of \cite[Proposition~5.2.1.17]{Lurie2},
we obtain the following corollary
by Theorems~\ref{theorem:L-R-perfect-pairing}
and \ref{thm:perfect-pairing-TS},
and Lemma~\ref{lemma:map-pairing-representable}.

\begin{corollary}\label{cor:tw-fun-oplax-eq-compatible-orig}
There is a natural equivalence of $(\infty,1)$-categories
\[ \twrsl{\scriptx}^{\rm op}
   \stackrel{\simeq}{\longrightarrow}
   \laxarrowsr{\scriptx} \]
which makes the following diagram commute
\[ \begin{array}{ccc}
    \twrsl{\scriptx}^{\rm op}
    &\stackrel{\simeq}{\longrightarrow}&
    \laxarrowsr{\scriptx}\\
    \mbox{$\scriptstyle (s,t)_{sL}$}
    \bigg\downarrow
    \phantom{\mbox{$\scriptstyle (s,t)_{sL}$}}& &
    \phantom{\mbox{$\scriptstyle (s,t)_{sR}$}}
    \bigg\downarrow
    \mbox{$\scriptstyle (s,t)_{sR}$}\\
    (u_1\scriptx^L)^{\rm op}\times u_1\scriptx&
    \stackrel{\simeq}{\longrightarrow}&
    u_1\scriptx^R\times u_1\scriptx.\\
\end{array}\]
\end{corollary}

By applying Corollary~\ref{cor:tw-fun-oplax-eq-compatible-orig}
to the $(\infty,2)$-category $\mathcal{X}^{\twoop}$,
we obtain the following corollary.

\begin{corollary}\label{cor:tw-fun-lax-eq-compatible}
There is a natural equivalence of $(\infty,1)$-categories
\[ \twlsr{\scriptx}
   \stackrel{\simeq}{\longrightarrow}
   \oplaxarrowsl{\scriptx} \]
which makes the following diagram commute
\[ \begin{array}{ccc}
    \twlsr{\scriptx}
    &\stackrel{\simeq}{\longrightarrow}&
    \oplaxarrowsl{\scriptx}\\
    \mbox{$\scriptstyle (s,t)_{sR}$}
    \bigg\downarrow
    \phantom{\mbox{$\scriptstyle (s,t)_{sR}$}}
    & &
    \phantom{\mbox{$\scriptstyle (s,t)_{sL}$}}
    \bigg\downarrow\mbox{$\scriptstyle (s,t)_{sL}$} \\
    (u_1\scriptx^R)^{\rm op}
    \times u_1\scriptx
    &\stackrel{\simeq}{\longrightarrow}&
    u_1\scriptx^L\times u_1\scriptx.\\
\end{array}\]
\end{corollary}

\section{The $\mathcal{O}$-monoidal structure on
$\scriptz(A,A)^{\otimes}_{\bullet}$
and convolution product}
\label{section:duoidal-convolution}

Let $\scriptz$ be an $\mathcal{O}$-monoidal
$(\infty,2)$-category
and let $A$ be a map $\mathcal{O}$-monoidale in $\scriptz$.
In \S\ref{section:map-duoidal}
we showed that
$\scriptz(A,A)^{\otimes}_{\bullet}$
corresponds to
a coCartesian $(\Delta^{\rm op},\mathcal{O})$-duoidal
$\infty$-category.
In this section
we show that the $\mathcal{O}$-monoidal
structure on
$\scriptz(A,A)^{\otimes}_1$
given by Theorem~\ref{thm:Z-A-A-duoidal-structure}
is equivalent to the convolution product
$\scriptz(A^*,A)^{\otimes}_{\rm conv}$
defined in \S\ref{section:convolution_product}.

First,
we promote the equivalences in
Corollaries~\ref{cor:tw-fun-oplax-eq-compatible-orig}
and \ref{cor:tw-fun-lax-eq-compatible}
to those of Cartesian symmetric monoidal $\infty$-categories.
We consider the functors
$\TSoprlpair(-)$,
$\twrsl{-}$,
$\laxarrowsr{-}$,
${\rm Pair}(-)\times \twl{(-)^{\twoop}}$,
and
$(u_1(-)^L\times u_1(-)^R)\times
(u_1(-)^{\rm op}\times u_1(-))$
from $\bicat$ to $\cat$.
Since these functors preserve finite products,
they extend to symmetric monoidal
functors
between
the Cartesian symmetric monoidal
$\infty$-categories $\bicat^{\times}$ and $\cat^{\times}$.
Furthermore,
the natural transformations
appearing in Corollary~\ref{cor:tw-fun-oplax-eq-compatible-orig}
extend to natural transformations
of symmetric monoidal functors.
Hence,
we obtain the following proposition.

\begin{proposition}\label{prop:monoidal-tw-arrow-eq}
For an $\mathcal{O}$-monoidal
$(\infty,2)$-category $\scriptz$,
there is a natural equivalence
\[ \twrsl{\scriptz}^{\otimes}
   \stackrel{\simeq}{\longrightarrow}
   \laxarrowsr{\scriptz}^{\otimes}  \]
of $\mathcal{O}$-monoidal $(\infty,1)$-categories,
which makes the following diagram commute
\[ \begin{array}{ccc}
    \twrsl{\scriptz}^{\otimes}
    &\stackrel{\simeq}{\longrightarrow}&
    \laxarrowsr{\scriptz}^{\otimes}\\
    \mbox{$\scriptstyle (s,t)_{sL}^{\otimes}$}
    \bigg\downarrow
    \phantom{\mbox{$\scriptstyle (s,t)_{sL}^{\otimes}$}}
    & &
    \phantom{\mbox{$\scriptstyle (s,t)_{sR}^{\otimes}$}}
    \bigg\downarrow
    \mbox{$\scriptstyle (s,t)_{sR}^{\otimes}$}\\
    ((u_1\scriptz^L)^{\rm op})^{\otimes}
    \times u_1\scriptz^{\otimes}&
    \stackrel{\simeq}{\longrightarrow}&
    (u_1\scriptz^R)^{\otimes}\times
    u_1\scriptz^{\otimes}.\\
\end{array}\]
\end{proposition}

By applying Proposition~\ref{prop:monoidal-tw-arrow-eq}
to the $2$-opposite $\scripto$-monoidal
$(\infty,2)$-category $(\scriptz^{\twoop})^{\otimes}$,
we obtain the following proposition.

\begin{proposition}\label{prop:monoidal-tw-fun-eq}
For an $\mathcal{O}$-monoidal
$(\infty,2)$-category $\scriptz$,
there is a natural equivalence
\[ \twlsr{\scriptz}^{\otimes}
   \stackrel{\simeq}{\longrightarrow}
   \oplaxarrowsl{\scriptz}^{\otimes}  \]
of $\mathcal{O}$-monoidal $(\infty,1)$-categories,
which makes the following diagram commute
\[ \begin{array}{ccc}
    \twlsr{\scriptz}^{\otimes}
    &\stackrel{\simeq}{\longrightarrow}&
    \oplaxarrowsl{\scriptz}^{\otimes}\\
    \mbox{$\scriptstyle (s,t)_{sR}^{\otimes}$}
    \bigg\downarrow
    \phantom{\mbox{$\scriptstyle (s,t)_{sR}^{\otimes}$}}
    & &
    \phantom{\mbox{$\scriptstyle (s,t)_L^{\otimes}$}}
    \bigg\downarrow
    \mbox{$\scriptstyle (s,t)_L^{\otimes}$}\\
    ((u_1\scriptz^R)^{\rm op})^{\otimes}
    \times u_1\scriptz^{\otimes}&
    \stackrel{\simeq}{\longrightarrow}&
    (u_1\scriptz^L)^{\otimes}\times
    u_1\scriptz^{\otimes}.\\
\end{array}\]
\end{proposition}

Next,
we consider the $\mathcal{O}$-monoidal
$(\infty,1)$-category $\scriptz(A,A)^{\otimes}_1$
for a map $\mathcal{O}$-monoidale $A$ in $\scriptz$.
We shall slightly generalize
the construction of $\scriptz(A,A)^{\otimes}_1$.
Let $B\in {\rm Alg}_{\mathcal{O}}(u_1\scriptz)$
be an $\mathcal{O}$-monoidale in $\scriptz$.
We define an $\infty$-operad
$\scriptz(A,B)^{\boxtimes}$
by the following pullback diagram
\[ \begin{array}{ccc}
  \scriptz(A,B)^{\boxtimes}
  & \longrightarrow &
  \oplaxarrowsl{\scriptz}^{\otimes}\\
  \bigg\downarrow & &
  \phantom{\mbox{$\scriptstyle (s,t)_{sL}^{\otimes}$}}
  \bigg\downarrow
  \mbox{$\scriptstyle (s,t)_{sL}^{\otimes}$}\\
  \scripto^{\otimes}
  & \stackrel{(A,B)}{\longrightarrow} &
  (u_1\scriptz^L)^{\otimes}\times
  (u_1\scriptz)^{\otimes}
   \end{array}\]
By Corollary~\ref{cor:Ar-sL-coCart},
$\scriptz(A,B)^{\boxtimes}$
is an $\scripto$-monoidal $\infty$-category.
We notice that
there is an equivalence
$\scriptz(A,A)^{\boxtimes}\simeq
\scriptz(A,A)^{\otimes}_1$
of $\scripto$-monoidal
$\infty$-categories.

By Proposition~\ref{prop:monoidal-tw-fun-eq},
we obtain the following theorem
which means that the $\mathcal{O}$-monoidal
structure on $\scriptz(A,B)^{\boxtimes}$
is equivalent to that given by convolution product.

\begin{theorem}
Let $A\in {\rm Alg}_{\mathcal{O}}(u_1\scriptz^{L})$
be a map $\mathcal{O}$-monoidale 
and let 
$B\in {\rm Alg}_{\mathcal{O}}(u_1\scriptz)$
be an $\mathcal{O}$-monoidale in $\scriptz$.
We denote by $A^*\in {\rm coAlg}_{\mathcal{O}}(u_1\scriptz^R)$
the right adjoint $\mathcal{O}$-comonoidale associated to $A$.
There is an equivalence of $\mathcal{O}$-monoidal $\infty$-categories
\[ \scriptz(A^*,B)^{\otimes}_{\rm conv}\simeq
   \scriptz(A,B)^{\boxtimes}. \]
\end{theorem}

\proof
There is a commutative diagram
\[ \xymatrix{
    & \mathcal{O}^{\otimes}\ar[dl]_{(A^*,B)}\ar[dr]^{(A,B)} & \\
  ((u_1\scriptz^R)^{\rm op})^{\otimes}\times u_1\scriptz^{\otimes}
  \ar[rr]^{\simeq} & &
  (u_1\scriptz^{L})^{\otimes}\times u_1\scriptz^{\otimes}
}\]
of $\infty$-operads.
The theorem follows  
by taking pullbacks along $(A^*,B)$ and $(A,B)$
in the commutative diagram in
Proposition~\ref{prop:monoidal-tw-fun-eq}.
\qed

\begin{corollary}
Let $A\in {\rm Alg}_{\mathcal{O}}(u_1\scriptz^{L})$
be a map $\mathcal{O}$-monoidale 
in $\scriptz$.
There is an equivalence
of $\mathcal{O}$-monoidal $\infty$-categories
\[ \scriptz(A^*,A)^{\otimes}_{\rm conv}\simeq
   \scriptz(A,A)^{\otimes}_1. \]
\end{corollary}

\section{Proof of Proposition~\ref{prop:TSZ-complete-Segal-space-orig}}
\label{section:proof-twisted-sq-complete-Segal}

In this section
we give a proof of
Proposition~\ref{prop:TSZ-complete-Segal-space-orig},
which says that
$\TSop(X)_{\bullet}$ is a complete Segal space
for any $\infty$-bicategory $X$.

First,
we show that
$\TSop(X)_{\bullet}$ is a Segal space.

\begin{proposition}
\label{prop:twisted-square-Segal-space}
For any $\infty$-bicategory $X$,
the simplicial object $\TSop(X)_{\bullet}$
is a Segal space.
\end{proposition}

In order to prove
Proposition~\ref{prop:twisted-square-Segal-space},
we recall the following lemma.

\begin{lemma}[{cf.~\cite[Lemma~4.6]{Torii3}}]
\label{lemma:general-reedy-Segal-space}
Let $A^{\bullet}$ be a cosimplicial object of $\setsc$,
and let $X$ be an $\infty$-bicategory.
If $A^{\bullet}$ is Reedy cofibrant
and satisfies the co-Segal condition,
then
the simplicial object
$\mapsc(A^{\bullet},X)$
is a Segal space.
\end{lemma}

By Lemma~\ref{lemma:general-reedy-Segal-space},
in order to prove
Proposition~\ref{prop:twisted-square-Segal-space},
it suffices to show that the cosimplicial
scaled simplicial set
$\tsscbullet$ is Reedy cofibrant
and satisfies the co-Segal condition.

\begin{lemma}\label{lemma:ts-sc-Reedy-cofibrant}
The cosimplicial scaled simplicial set 
$\tsscbullet$
is Reedy cofibrant.
\end{lemma}

\proof
The lemma follows by observing that
the underlying simplicial set of
the $n$th latching object of $\tsbullet$
is isomorphic to 
$(\Delta^2\times\partial\Delta^n)\cup
\Omega(\Delta^2\times\partial \Delta^n)$,
which is a subcomplex of
the underlying simplicial set
$(\Delta^2\times\Delta^n)\cup
\Omega(\Delta^2\times\Delta^n)$ of $\tsscn$.
\qed

\bigskip

Next, we would like to show that 
$\tsbullet$ satisfies
the co-Segal condition.
For this purpose,
we make some preliminary results.
For $0<i<n$,
we set
\[ \Lambda^n_i \tssc
   =((\Delta^2\times\Lambda^n_i)
   \cup
   \Omega(\Delta^2\times \Lambda^n_i))_{\dagger}.
   \]
We show that
the inclusion map
$\Lambda^n_i \tssc\to \tsscn$
is scaled anodyne.
In order to prove this,
we recall the following lemma,
which is a special case of
\cite[Lemma~1.18]{Garcia-Stern}.
Before that,
we recall the notation.
For a nonempty finite totally ordered set $S$
and a subset $\emptyset \neq N\varsubsetneqq S$,
we set
\[ \Lambda^S_N=\bigcup_{s\in S-N}\Delta^{S-\{s\}}. \]

\begin{lemma}[{\cite[Lemma~3.7]{Torii3},
\cite[Lemma~1.18]{Garcia-Stern}}]
\label{lemma:simplex-scaled-anodyne}
Let $\Delta^r_{\dagger}=(\Delta^r,T)$
be a scaled simplicial set for $r\ge 3$,
and let $M$ be a nonempty subset of
$\{0,1,\ldots,r-1\}$.
We assume that
there exists an integer $s$ with
$0\le s<t$
such that
$s\not\in M$ and $i\in M$ for all $s<i\le t$,
where $t$ is the largest number of $M$.
Furthermore,
we assume
that $|M|\le r-2$ and
that the triangle $\Delta^{[r]-M}$
is not thin in $\Delta^r_{\dag}$
when $|M|=r-2$.
If triangles $\Delta^{\{i,t,t+1\}}$
are thin in $\Delta^r_{\dag}$ for all $s\le i< t$,
then the inclusion map
$\Lambda^r_{M,\dag}\to \Delta^r_{\dag}$
is scaled anodyne,
where $\Lambda^r_{M,\dagger}$ is $\Lambda^r_M$
equipped with the induced scaling from $\Delta^r_{\dagger}$.
\end{lemma}

We decompose $\Lambda^n_i\tssc$
into two parts.
We set 
\[ \begin{array}{rcl}
    \Lambda^n_i \tsscplusno
    &=&
    (\Delta^2\times\Lambda^n_i)_{\dagger},\\[2mm]
    \widehat{\Lambda}^n_i\tsscminusno
    &=&
    \Omega((\Delta^2\times\Lambda^n_i)\cup
    (\Delta^{\{0,2\}}\times \Delta^n))_{\dagger}.\\
   \end{array}   \]

\begin{lemma}\label{lemma:scaled-anodyne-lambda-plus}
For any $0<i<n$,
the inclusion map
$\Lambda^n_i \tsscplusno
\to\tsscnplus$
is scaled anodyne.
\end{lemma}

\proof
First,
we shall show that the inclusion map
$\Lambda^n_i\tsscplusno
   \to
   \overline{\Lambda}^n_i \tsscplusno$
is scaled anodyne,
where
$\overline{\Lambda}^n_i
   \tsscplusno=
   \Lambda^n_i\tsscplusno
   \cup (\Lambda^2_1\times \Delta^n)_{\dagger}$.
Since
$\Lambda^n_{i,\sharp}\to\Delta^n_{\sharp}$
is a scaled anodyne map of
type (An1),
the inclusion map
$\Lambda^n_i\tsscplusno\to
\Lambda^n_i\tsscplusno
\cup \cup_{i=0,1,2}(\Delta^{\{i\}}\times \Delta^n)_{\dagger}$
is scaled anodyne.
By using
\cite[Proposition~2.16]{GHL2} and 
isomorphisms
$(\Delta^{\{0,1\}}\times\Delta^n)_{\dagger}
\cong \Delta^1_{\sharp}\otimes \Delta^n_{\sharp}$
and
$(\Delta^{\{1,2\}}\times \Delta^n)_{\dagger}\cong
\Delta^n_{\sharp}{\otimes} \Delta^1_{\sharp}$,
we see that
the inclusion map
$\Lambda^n_i\tsscplusno
\cup \cup_{i=0,1,2}(\Delta^{\{i\}}\times \Delta^n)_{\dagger}
\to
\overline{\Lambda}^n_i\tsscplusno$
is scaled anodyne.

Next,
we shall show that the inclusion map
$\overline{\Lambda}^n_i
   \tsscplusno \to \tsscnplus$
is scaled anodyne.
For $0\le r\le n$ and $0\le k\le n-r$,
we denote by $\sigma(r,k)$
the $(n+2)$-dimensional simplex
\[\Delta^{\{000,\ldots,00k,01k,\ldots,01(k+r),11(k+r),\ldots,11n\}}\]
of $\Delta^2\times \Delta^n$.
For $0\le s\le n$,
we set
\[ X(s)= \overline{\Lambda}^n_i \tsscplusno
         \cup  \bigcup_{s\le r\le n}\
         \bigcup_{0\le k\le n-r} \sigma(r,k)_{\dagger}.\]
Then,
we obtain a filtration
\[ \overline{\Lambda}^n_i \tsscplusno=
   X(n+1)\to X(n)\to
   X(n-1)\to\cdots
   \to X(0)
   = \tsscnplus .\]
It suffices to show that the inclusion map
$X(s+1)\to X(s)$
is scaled anodyne for each $0\le s\le n$.

We fix $s$ with $0\le s\le n$
and show that $X(s+1)\to X(s)$ is scaled anodyne.
We identify $\sigma(s,k)_{\dagger}$
with a map $\sigma(s,k): \Delta^{n+2}_{\dagger}\to X(s)$,
where $\Delta^{n+2}_{\dagger}$ is $\Delta^{n+2}$
equipped with the induced scaling.
There is a pushout diagram
\[ \begin{array}{ccc}
  \coprod_{0\le k\le n-s}\Lambda^{n+2}_{M(k),\dagger}
  & \longrightarrow &
     \coprod_{0\le k\le n-s}\Delta^{n+2}_{\dagger}\\
     \bigg\downarrow & &
     \phantom{\mbox{$\scriptstyle \sigma(s,k)$}}
     \bigg\downarrow
     \mbox{$\scriptstyle \sigma(s,k)$}\\
     X(s+1) & \longrightarrow & X(s),
   \end{array}\]  
where
\[ M(k)=\left\{\begin{array}{ll}
                \{00i,01k,01(k+s)\} & (0<i<k),\\
                \{01k,01i,01(k+s)\} & (k\le i\le k+s),\\
                \{01k,01(k+s),11i\} & (k+s<i<n).\\
               \end{array} \right. \]
Thus,
it suffices to show that
the inclusion map
$\Lambda^{n+2}_{M(k),\dagger}\to \Delta^{n+2}_{\dagger}$
is scaled anodyne. 
This follows from Lemma~\ref{lemma:simplex-scaled-anodyne}.
\qed

\bigskip

\begin{lemma}\label{lemma:scaled-anodyne-minus}
For any $0<i<n$,
the inclusion map
$\widehat{\Lambda}^n_i
   \tsscminusno
   \to
   \tsscnminus$
is scaled anodyne.
\end{lemma}

\proof
First,
we shall show that
the inclusion map
$\widehat{\Lambda}^n_i\tsscminusno
   \to
   \overline{\Lambda}^n_i\tsscminusno$
is scaled anodyne,
where
$\overline{\Lambda}^n_i \tsscminusno
    =    
 \widehat{\Lambda}^n_i \tsscminusno 
    \cup \Omega(\Lambda^2_1\times \Delta^n)_{\dagger}$.
Since
$\Lambda^n_{i,\sharp}\to\Delta^n_{\sharp}$
is a scaled anodyne map of
type (An1),
the inclusion map
$\widehat{\Lambda}^n_i\tsscminusno\to
\widehat{\Lambda}^n_i\tsscminusno
\cup \Omega(\Delta^{\{1\}}\times \Delta^n)_{\dagger}$
is scaled anodyne.
By the proof of \cite[Lemma~5.5]{Torii3} and its dual,
we see that
the inclusion map
$\widehat{\Lambda}^n_i\tsscminusno
\cup \Omega(\Delta^{\{1\}}\times \Delta^n)_{\dagger}
\to
\overline{\Lambda}^n_i\tsscminusno$
is scaled anodyne.

Next,
we shall show that
the inclusion map
$\overline{\Lambda}^n_i\tsscminusno
\to \tsscnminus$
is scaled anodyne.
For $0\le r\le n$ and $0\le k\le n-r$,
we denote by $\overline{\sigma}(r,k)$
the $(n+2)$-dimensional simplex
\[\Delta^{\{000,\ldots,00k,10k,\ldots,10(k+r),11(k+r),\ldots,11n\}}\]
of $\Delta^n\star\Delta^{n,{\rm op}}\star\Delta^n$.
For $0\le s\le n$,
we set
\[ Y(s)= \overline{\Lambda}^n_i\tsscminusno
         \cup  \bigcup_{0\le r\le s}\
         \bigcup_{0\le k\le n-r} \overline{\sigma}(r,k)_{\dagger}.\]
Then,
we obtain a filtration
\[ \overline{\Lambda}^n_i\tsscminusno=
   Y(-1)\to Y(0)\to
   Y(1)\to\cdots
   \to Y(n)
   = \tsscnminus.\]
It suffices to show that the inclusion map
$Y(s-1)\to Y(s)$
is scaled anodyne for each $0\le s\le n$.

We fix $s$ with $0\le s\le n$
and show that $Y(s-1)\to Y(s)$ is scaled anodyne.
We identify $\overline{\sigma}(s,k)_{\dagger}$
with a map $\overline{\sigma}(s,k):
\Delta^{n+2}_{\dagger}\to Y(s)$,
where $\Delta^{n+2}_{\dagger}$ is $\Delta^{n+2}$
equipped with the induced scaling.
There is a pushout diagram
\[ \begin{array}{ccc}
  \coprod_{0\le k\le n-s}\Lambda^{n+2}_{M(k)}
  & \longrightarrow &
     \coprod_{0\le k\le n-s}\Delta^{n+2}_{\dagger}\\
     \bigg\downarrow & &
     \phantom{\mbox{$\scriptstyle \sigma(s,k)$}}
     \bigg\downarrow
     \mbox{$\scriptstyle \overline{\sigma}(s,k)$}\\
     Y(s-1) & \longrightarrow & Y(s),
   \end{array}\]  
where
\[ M(k)=\left\{\begin{array}{ll}
                \{10i,11s\} & (0=k< i<s<n),\\
                \{10i\} & (0=k< i<s=n),\\
                \{11s,11i\} & (0=k\le s\le i<n),\\
                \{00i,00k,11(k+s)\} & (0<i\le k\le k+s<n),\\
                \{00k,10i,11(k+s)\} & (0<k< i< k+s<n),\\
                \{00k,11(k+s),11i\} & (0<k\le k+s\le i<n),\\
                \{00i,00k\} & (0<i\le k\le k+s=n),\\
                \{00k,10i\} & (0<k< i< k+s=n).\\

\end{array} \right. \]
Thus,
it suffices to show that
the inclusion map
$\Lambda^{n+2}_{M(k),\sharp}\to \Delta^{n+2}_{\dagger}$
is scaled anodyne. 
This follows from Lemma~\ref{lemma:simplex-scaled-anodyne}.
\qed

\begin{lemma}\label{lemma:inner-tssc}
For $0< i< n$,
the inclusion map
$\Lambda^n_i\tssc \to \tsscn$
is scaled anodyne.
\end{lemma}

\proof
The map
$\Lambda^n_i\tssc\to \tsscn$
can be written as the composite
$\Lambda^n_i\tssc\to
\Lambda^n_i\tssc\cup \tsscnplus
\to \tsscn$.
The first and second maps are scaled anodyne
since they are obtained as pushouts
of the maps in 
Lemmas~\ref{lemma:scaled-anodyne-lambda-plus}
and \ref{lemma:scaled-anodyne-minus},
respectively.
\qed

\proof
[Proof of Proposition~\ref{prop:twisted-square-Segal-space}]
By Lemmas~\ref{lemma:general-reedy-Segal-space}
and \ref{lemma:ts-sc-Reedy-cofibrant},
it suffices to show that
the cosimplicial object
$\tsscbullet$
satisfies the co-Segal condition.
We can prove 
that the co-Segal map
$\tsscone 
   \coprod_{\tssczero}\cdots
   \coprod_{\tssczero}
   \tsscone
   \to\tsscn$
is scaled anodyne 
by induction on $n$
together with Lemma~\ref{lemma:inner-tssc}.
\qed

\bigskip

Next,
we show that
the Segal space
$\TSop(X)_{\bullet}$ is complete.
In order to prove this,
we make some preliminary results.

By Definition~\ref{def:partial-maps},
we have maps of Segal spaces
\[ \partial^f: \TSop(X)_{\bullet}\longrightarrow
   \partial^f\TSop(X)_{\bullet} \]
for $f={\rm T}, {\rm F}, {\rm R}, {\rm B}$. 
We set
\[ \bd(X)_{\bullet}=
   \prod_f\partial^f\TSop(X)_{\bullet},\]
where the product is taken over
$f={\rm T}, {\rm F}, {\rm R}, {\rm B}$. 
Note that $\bd(X)_{\bullet}$ is a complete Segal space,
that is,
the map $s_0: \bd(X)_1\to \bd(X)_1^{\rm eq}$
is an equivalence,
where $\bd(X)_1^{\rm eq}$ is the full subspace
of $\bd(X)_1$ spanned by equivalences. 

In order to prove that
the Segal space $\TSop(X)_{\bullet}$
is complete,
it suffices to show that 
the map
$s_0: \TSop(X)_0 \to \TSop(X)^{\rm eq}_1$
is an equivalence.
We denote by 
$\partial^{\rm Bd}: \TSop(X)_{\bullet}\to \bd(X)_{\bullet}$
the map of Segal spaces
obtained as a product of $\partial^f$
for $f={\rm T}$, ${\rm F}$, ${\rm R}$, ${\rm B}$.
We define
$\TSop(X)_1^{{\rm bd},{\rm eq}}$
to be 
\[ \TSop(X)_1\times_{\bd(X)_1} \bd(X)_1^{\rm eq}. \]
We notice that
the map $s_0: \TSop(X)_0\to \TSop(X)_1$
factors through $\TSop(X)_1^{\rm bd,eq}$.

\begin{lemma}
\label{lemma:s0-eq-if-s0-bd-eq}
If $s_0: \TSop(X)_0\to \TSop(X)_1^{\rm bd,eq}$
is an equivalence,
then $s_0: \TSop(X)\to \TSop(X)^{\rm eq}$
is also an equivalence.
\end{lemma}

\proof
The lemma follows from
the fact that
$\TSop(X)^{\rm eq}_1$ is a full subspace
of $\TSop(X)^{\rm bd,eq}_1$.
\qed

\bigskip

We consider the following commutative diagram
\begin{equation}\label{eq:complete-fiber-sqaure-diagram}
\begin{array}{ccc}
    \TSop(X)_0
    &\stackrel{s_0}{\longrightarrow}&
    \TSop(X)^{{\rm bd},{\rm eq}}_1 \\
    \bigg\downarrow & & \bigg\downarrow \\
    \bd(X)_0 & \stackrel{s_0}{\longrightarrow}&
    \bd(X)_1^{\rm eq}.\\
   \end{array}
\end{equation}
Since the bottom horizontal arrow is an equivalence,
it suffices to show that
diagram (\ref{eq:complete-fiber-sqaure-diagram})
is pullback
in order to prove that the top horizontal arrow is
an equivalence.
For $P\in \bd(X)_0$,
let
$(s_0)_P: \TSop(X)_{0,P}\to
   \TSop(X)_{1,s_0(P)}^{{\rm bd},{\rm eq}}$
be the induced map on fibers.
We show that
$(s_0)_P$ is an equivalence
for any $P\in \bd(X)_0$.
In fact,
we show that
$(d_i)_P: \TSop(X)_{1,s_0(P)}^{{\rm bd},{\rm eq}}
\to \TSop(X)_{0,P}$
is an equivalence
for $i=0,1$.

We consider
a scaled simplicial set
\[ \wtsscone=
   ((\Delta^2\times\Delta^1)\cup
     \Omega(\Delta^2\times\Delta^1),
     \widetilde{T}),\]
where
$\widetilde{T}$ is a set of $2$-simplices
that consists of the thin triangles of $\tsscone$
together with
\[ \begin{array}{lll}
   \Delta^{\{000,001,011\}},&
   \Delta^{\{010,110,111\}},&
   \Delta^{\{000,100,101\}},\\
   \Delta^{\{100,101,111\}},&
   \Delta^{\{000,001,111\}},&
   \Delta^{\{000,110,111\}}.\\
   \end{array}\]
Note that there is an inclusion map
$\tsscone\to \wtsscone$
of scaled simplicial sets.

For a subcomplex $K$ of
(the underlying simplicial set of) $\tsscone$,
we write $K_{\diamondsuit}$ for
the scaled simplicial set $K$
equipped with the induced scaling
from $\wtsscone$.
We set 
\[ {\rm{FSR}}^i=
   (\partial^{\rm F}\tsscone\cup
   d^i(\tssczero) \cup
   \partial^{\rm R}\tsscone)_{\diamondsuit}\]
for $i=0,1$.
We consider inclusion maps
\[ \begin{array}{lrcl}
   \theta_0:&
   {\rm FSR}^0   \coprod_{\Delta^{\{110,111\}}_{\sharp}}
   \Delta^{\{0\}}_{\sharp}
   &\longrightarrow&
   \wtsscone\coprod_{\Delta^{\{110,111\}}_{\sharp}}
   \Delta^{\{0\}}_{\sharp},\\
    \theta_1:&
   {\rm FSR}^1
   \coprod_{\Delta^{\{000,001\}}_{\sharp}}\Delta^{\{0\}}_{\sharp}
   &\longrightarrow& 
   \wtsscone\coprod_{\Delta^{\{000,001\}}_{\sharp}}
   \Delta^{\{0\}}_{\sharp}.\\
   \end{array}\]

\begin{lemma}\label{lemma:theta-trivial-cof}
The inclusion maps
$\theta_0$ and $\theta_1$
are trivial cofibrations in $\setsc$.
\end{lemma}

\proof
We shall show that the map
$\theta_1$
is a trivial cofibration.
The other case is similar.
For simplicity,
we set
\[ \begin{array}{rcl}
    E_0&=&{\rm FSR}^1\coprod_{\Delta^{\{000,001\}}_{\sharp}}
    \Delta^{\{0\}}_{\sharp},\\
     E_1&=&({\rm FSR}^1\cup
      \tssconeminus{}_{\diamondsuit})
     \coprod_{\Delta^{\{000,001\}}_{\sharp}}
     \Delta^{\{0\}}_{\sharp},\\
    E_2&=&\tsscone_{\diamondsuit}
    \coprod_{\Delta^{\{000,001\}}_{\sharp}}
    \Delta^{\{0\}}_{\sharp}.\\
  \end{array}\]    
It suffices to show that
the inclusion maps
$E_0\to E_1$ and $E_1\to E_2$ are trivial cofibrations.

First,
we shall show that
$E_0\to E_1$ is a trivial cofibration.
We set
\[ \begin{array}{rcl}
  F_0&=& {\rm FSR}^1,\\
  F_1&=& {\rm FSR}^1\cup
         \Delta^{\{000,100,110,111\}}_{\diamondsuit} ,\\
  F_2&=& {\rm FSR}^1\cup (\Delta^{\{000,100,110,111\}}\cup
  \Delta^{\{000,100,101,111\}})_{\diamondsuit} .\\
\end{array}\]
We can see that the inclusion maps
$F_0\to F_1$ and $F_1\to F_2$ 
are scaled anodyne
since they are obtained by iterated
pushouts along scaled anodyne maps
of type (An1).
By taking a pushout along the map
$\Delta^{\{000,001\}}_{\sharp}\to \Delta^{\{0\}}_{\sharp}$,
we see that
$E_0\to F_2\coprod_{\Delta^{\{000,001\}}_{\sharp}}\Delta^0_{\sharp}$
is scaled anodyne.
By \cite[Remark~1.33]{GHL},
the map
$\Lambda^2_{0,\sharp}\coprod_{\Delta^{\{0,1\}}_{\sharp}}
 \Delta^{\{0\}}_{\sharp}\to
 \Delta^2_{\sharp}\coprod_{\Delta^{\{0,1\}}_{\sharp}}
 \Delta^{\{0\}}_{\sharp}$
is a trivial cofibration.
By taking a pushout along it,
we obtain a trivial cofibration
$F_2\coprod_{\Delta^{\{000,001\}}_{\sharp}}\Delta^0_{\sharp}
\to E_1$.
Hence,
$E_0\to E_1$ is a trivial cofibration
by composition.

Next,
we shall show that $E_1\to E_2$ is a trivial cofibration.
We set
\[ \begin{array}{rcl}
  G_0&=& {\rm FSR}^1\cup
  \tssconeminus{}_{\diamondsuit},\\
  G_1&=& {\rm FSR}^1\cup
  \tssconeminus{}_{\diamondsuit}\cup
  \Delta^{\{000,010,110,111\}}_{\diamondsuit},\\
  G_2&=& {\rm FSR}^1\cup
  \tssconeminus{}_{\diamondsuit}\cup
  (\Delta^{\{000,010,110,111\}}\cup
  \Delta^{\{000,010,011,111\}})_{\diamondsuit}.\\
   \end{array} \]
We can see that the inclusion maps
$G_0\to G_1$ and $G_1\to G_2$
are scaled anodyne
since they are obtained by taking pushouts
along scaled anodyne maps of type (An1).
By taking a pushout along the map
$\Delta^{\{000,001\}}_{\sharp}\to \Delta^{\{0\}}_{\sharp}$,
we see that
$E_1\to G_2\coprod_{\Delta^{\{000,001\}}_{\sharp}}\Delta^0_{\sharp}$
is scaled anodyne.
By taking a pushout along
the trivial cofibration
$\Lambda^2_{0,\sharp}\coprod_{\Delta^{\{0,1\}}_{\sharp}}
 \Delta^{\{0\}}_{\sharp}\to
 \Delta^2_{\sharp}\coprod_{\Delta^{\{0,1\}}_{\sharp}}
 \Delta^{\{0\}}_{\sharp}$,
we obtain a trivial cofibration
$G_2\coprod_{\Delta^{\{000,001\}}_{\sharp}}\Delta^0_{\sharp}
\to E_2$.
Hence,
$E_1\to E_2$ is a trivial cofibration
by composition.
\qed

\bigskip

For $i=0,1$,
Lemma~\ref{lemma:theta-trivial-cof}
implies that the induced map 
\[ \Delta^1_{\sharp}\coprod_{\partial^{\rm F}\wtsscone}
   {\rm FSR}^i
    \coprod_{\partial^{\rm R}\wtsscone}\Delta^1_{\sharp}
    \longrightarrow
    \Delta^1_{\sharp}\coprod_{\partial^{\rm F}\wtsscone}
    \wtsscone   
    \coprod_{\partial^{\rm R}\wtsscone}\Delta^1_{\sharp} \]
is a trivial cofibration,
where $\partial^{\rm F}\tsscone_{\diamondsuit}\to \Delta^1_{\sharp}$
is a map given by
$010,011\mapsto 0$
and
$110,111\mapsto 1$,
and $\partial^{\rm R}\tsscone_{\diamondsuit}\to \Delta^1_{\sharp}$
is a map given by
$000,001\mapsto 0$
and
$100,101\mapsto 1$.
Since the map
$d^i: \tssczero\to \tsscone$
induces an isomorphism
of scaled simplicial sets
between
$\tssczero$ and 
$\Delta^1_{\sharp}\coprod_{\partial^{\rm F}\wtsscone}
   {\rm FSR}^i
   \coprod_{\partial^{\rm R}\wtsscone}\Delta^1_{\sharp}$,
it induces
a trivial cofibration
\[ \tssczero
   \longrightarrow 
   \Delta^1_{\sharp}\coprod_{\partial^{\rm F}\wtsscone}
   \wtsscone
  \coprod_{\partial^{\rm R}\wtsscone}\Delta^1_{\sharp} \]
in $\setsc$.
This implies that the restriction map
$d_i:    \TSop(X)_{1,(s_0(b_{000}),s_0(b_{010}))}^{\rm eq}
   \to
\TSop(X)_{0,(b_{000},b_{010})}$
is an equivalence.
Hence,
we obtain that
$(d_i)_P:  \TSop(X)_{1,s_0(P)}^{\rm bd, eq}
\to\TSop(X)_{0,P}$
is an equivalence.
Since the composite $(d_i)_p(s_0)_P$ is an identity,
we obtain the following corollary.

\begin{corollary}\label{cor:the-fiber-map-s0-equivalence}
The map
$(s_0)_P: \TSop(X)_{0,P}\to
   \TSop(X)_{1,s_0(P)}^{{\rm bd},{\rm eq}}$
is an equivalence
for any $P\in\bd(X)_0$.
\end{corollary}

\proof[Proof of
Proposition~\ref{prop:TSZ-complete-Segal-space-orig}]
By Proposition~\ref{prop:twisted-square-Segal-space}
and Lemma~\ref{lemma:s0-eq-if-s0-bd-eq},
it suffices to show that
the map $s_0: \TSop(X)_0\to \TSop(X)_1^{\rm bd,eq}$
is an equivalence.
By Corollary~\ref{cor:the-fiber-map-s0-equivalence},
(\ref{eq:complete-fiber-sqaure-diagram})
is a pullback diagram.
Since the bottom horizontal arrow
in (\ref{eq:complete-fiber-sqaure-diagram})
is an equivalence,
the top horizontal arrow is also an equivalence.
\qed

\end{document}